\documentclass [a4paper, 10pt]{amsart}
\usepackage[textheight=25 cm, textwidth=16cm, headheight=10pt, headsep=20pt, footskip=10pt]{geometry}

\usepackage{amssymb,mathrsfs}
\usepackage{MnSymbol}
  
\newcommand{\cal}{\mathcal }

\newcommand{\R}{{\mathbb R}}

\renewcommand{\a}{\alpha}

\newcommand{\pa}{\partial}

\newcommand{\D}{\Delta}

\def\s{\sigma}
\def\r{{\rho}}
\def\n{\nu}
\def\p{{\pi}}

\def\AA{{\cal A}}

\def\AA{{\mathfrak{A}}}

\def\be{\begin{equation}}
\def\ee{\end{equation}}
\def\bea{\begin{eqnarray}}
\def\eea{\end{eqnarray}}

\def\nn{\nonumber}

\def\o{\omega}
\def\b{\beta}
\def\t{\tau}

\def\bC {\mathbf{C}}

\def\bK {\mathbf{K}}

\def\bQ {\mathbf{Q}}
\def\bR {\mathbf{R}}
\def\bS {\mathbf{S}}

\def\a {{\alpha}}
\def\b {{\beta}}

\def\l {{\lambda}}

\newcommand{\ba}{\begin{aligned}}
\newcommand{\ea}{\end{aligned}}

\newtheorem{Thm}{Theorem}[section]
\newtheorem{Rmk}[Thm]{Remark}
\newtheorem{Prop}[Thm]{Proposition}
\newtheorem{Cor}[Thm]{Corollary}

\newtheorem{Def}[Thm]{Definition}

\usepackage{color}

\newcommand{\norm}[1]{\lVert #1 \rVert}

\newcommand{\tr}{\mbox{Tr}}
\newcommand{\FF}{F}
\newcommand{\unN}{1,\dots, N}
\newcommand{\un}[1]{1,\dots, #1}
\newcommand{\suml}{\sum\limits}
\newcommand{\prodl}{\prod\limits}

\newcommand{\sumls}[1]{\suml_{\substack{#1}}}
\newcommand{\ej}{\mathcal E_j}
\newcommand{\uk}{\underline k}

\newcommand{\bcr}{}
\newcommand{\bcb}{}
\newcommand{\ec}{}
\newcommand{\bcrr}[1]{#1}

\begin{document}

\title[Size of chaos for 
mean field models]{\LARGE On the size of chaos  in the mean field dynamics
\vspace{0.8cm}}

\author[T. Paul]{\large Thierry Paul}
\address[T.P.]{CMLS, Ecole polytechnique, CNRS, Universit\'e Paris-Saclay, 91128 Palaiseau Cedex, France}
\email{thierry.paul@polytechnique.edu}

\author[M. Pulvirenti]{\large Mario Pulvirenti}
\address[M.P.]{International Research Center on the Mathematics and Mechanics of Complex Systems \\�MeMoCS, University of L'Aquila, Italy}  
\email{pulviren@mat.uniroma1.it}

\author[S. Simonella]{\large Sergio Simonella}
\address[S.S.]{ Zentrum Mathematik, TU M\"{u}nchen \\ Boltzmannstrasse 3, 85748 Garching -- Germany}
\email{s.simonella@tum.de}

\maketitle

\begin{abstract}
 We consider the error arising from the approximation of an $N$-particle dynamics with its description in terms of a one-particle kinetic equation. We estimate
the distance between the $j$-marginal of the system and the factorized state, obtained in a mean field limit as $N \to \infty$. Our analysis relies on the evolution equation for the ``correlation error"  rather than on the usual BBGKY hierarchy.  The rate of convergence is shown to be $O(j^2/ N)$ in any bounded interval of time (size of chaos), as expected from heuristic arguments.
Our formalism applies to an abstract hierarchical 
mean field model with bounded collision operator and a large class of initial data, covering (a) stochastic jump processes 
converging to the homogeneous Boltzmann and the Povzner equation and (b)  quantum systems giving rise to the Hartree equation.
\end{abstract}

\tableofcontents
\LARGE
\section{Introduction}

The kinetic description of particle systems is based on the propagation of chaos. This property allows to
substitute the complex dynamics of a huge number of particles by a single nonlinear partial differential equation
for the 
probability density (in quantum systems, the reduced density matrix) of a given particle. More precisely, one applies a statistical description.  At time zero,
the $N-$particle system is assumed to be ``chaotic" in the sense that each particle is distributed identically and independently from the others,  at least up to an error, vanishing when $N$ diverges. The dynamics creates correlations and the statistical independence is lost at any positive time. However,
after suitable rescaling of space and time, the statistical independence of any finite group of particles can be recovered, in the limit  $N \to \infty$. As a consequence a given particle 
evolves according to  an effective equation. The nature of this dynamics is
determined by the microscopic details of the system and by the regime of physical parameters. Such a mechanism 
works in the formal (in a few cases, rigorous) derivation of the most common kinetic equations. 

In this paper we consider the following class of mean field models.

\medskip
\noindent
$\bullet$ {\it Kac model.} The $N$-particle system evolves according to a stochastic process. To  each particle, say particle $i$, we associate a velocity $v_i \in \R^3$.
\newcommand{\calv}{\cal V}
The vector $\calv_N=\{v_1, \cdots, v_N\}$ changes by means of two-body collisions at random times, with random scattering angle.
The probability density $f^N(\calv_N,t)$ evolves according to the master equation (forward Kolmogorov equation) 
$$
\pa_t f^N=\frac 1N \sum_{i<j} \int d\omega B(\omega; v_i-v_j)\{f^N(\calv_N^{i,j})-f^N(\calv_N)\}\;,
$$
where $\calv_N^{i,j}=\{v_1, \cdots,v_{i-1}, v_i', v_{i+1}, \cdots, v_{j-1},v_j', v_{j+1}, \cdots, v_N\}$ and the pair $v_i' ,v_j'$ gives the outgoing velocities after a collision with scattering (unit) vector $\o$ and incoming velocities $v_i, v_j $. $\tfrac{B(\omega; v_i-v_j)}{|v_i-v_j|}$ is the differential cross-section of the two-body process.
The resulting kinetic equation is the homogeneous Boltzmann equation
$$
\pa_t f(v) = \int dv_1 \int d\omega B(\omega; v-v_1) \{ f(v') f(v_1') -f(v) f(v_1) \}\;.
$$
Such a model has been introduced by Kac \cite{Kac,Kac2} and has been largely investigated over recent times, see e.g.\,\cite {MM}.
Very similar stochastic systems including space variables and (space-)delocalized collisions are frequently used to justify numerical schemes 
\cite{PWZ,RW}. We will not mention them explicitly although they could be included in our analysis.

\smallskip
\noindent
$\bullet$  {\it `Soft spheres' model.} A slightly more realistic variant, taking into account
the positions of particles $X_N= \{x_1, \cdots, x_N\} \in \R^{3N}$ and relative transport,
was introduced by Cercignani \cite{Cer} and further investigated in \cite{LP}.
The  probability density $f^N(X_N,  \calv_N,t)$ evolves according to the equation
\bea
\pa_t f^N + \suml_{i=1}^N v_i \cdot \nabla_{x_i} f^N &=&
\frac 1N \suml_{i<j} \, h\left(|x_i-x_j|\right) B\left(\frac{x_i-x_j}{|x_i-x_j|}; v_i-v_j\right)
\nn \\
&&
\times\{f^N(X_N , \calv_N^{i,j})-f^N(X_N,\calv_N)\}\;. \nn
\eea
Here $h:\R^+ \to \R^+$ is a positive function with compact support. 
Now a pair of particles collides at a random distance with rate 
modulated by $h$. 
The associated kinetic equation is the Povzner equation
\bea
\pa_t f(x,v) + v\cdot \nabla_{x}f (x,v) &=&\int dv_1 \int dx_1\, h(|x-x_1|) B\left(\frac{x-x_1}{|x-x_1|}; v-v_1\right)  \nn \\
&&\times\{ f(x, v') f(x_1, v_1') -f(x, v) f(x_1, v_1) \}, \nn 
\eea
which can be seen as an $h-$mollification of the inhomogeneous Boltzmann equation (formally obtained when $h$
converges to a Dirac mass at the origin). 

\smallskip
\noindent
$\bullet$ {\it Quantum mean field model.} The $N$-particle quantum system has a mean field Hamiltonian
$$
H_N=- \frac {\hbar^2}2 \sum_{i=1}^N \D_{x_i}  + \frac 1N \sum_{i<j}�V(x_i-x_j),
$$
where $V$ is a  two-body potential on $\bR^d$ and $d$ is the dimension of the physical space.
A state of the system is a density matrix $\r^N$ whose time evolution is given by the von Neumann equation
$$
 \pa_t \r^N = \frac1{i\hbar}[H_N, \r^N ]\;
$$
(equivalent, modulo a global phase, to the   Schr\"odinger equation).
The effective equation for the one-particle density matrix $\r$ is the Hartree equation
$$
\pa_t \r = \frac1{i\hbar}[H_{MF}, \r ]
$$
where 
$$
H_{MF}=-\frac {\hbar^2}2 \D_{x} + \int dy V(x-y) \r (y,y)\;.
$$
Here $\r(x,y)$ is the kernel of $\r$ and hence $\r (y,y)$ is the spatial density.

\medskip

In all the above cases we assume symmetry in the particle labels. Moreover, we assume that
the initial state factorizes (or at least does so in the limit $N \to \infty$), namely $f^N(0)=f_0^{\otimes N}$ and
$\r^N(0)=\r_0^{\otimes N}$. At time $t>0$, in spite of the correlations created by the dynamics, 
the system  is still factorizing in the limit $N \to \infty$ through its $j$-particle marginals
$f^N_j$,  $\r^N_j$,  defined as partial integrations  of $f^N$ and partial  traces of $\r^N$, in the sense that these marginals converge, for any fixed $j$ and  in the limit $N \to \infty$,  to  $f^{\otimes j}$ and $\r^{\otimes j}$ respectively, 
$f=f(t)$ and $\r=\r(t)$ being the solutions of the associated effective equations.

This `propagation of chaos' has been  proved for the models under consideration and, under certain
assumptions, informations on the convergence rate are also available (see Section \ref{applications} below for 
bibliographical references).
\vskip 1cm

A natural question  arises ({\it size of chaos}): how large can be a group of $j=j(N)$ distinct particles, $j(N)$ diverging with $N$, so that one sees the decorrelation  of such systems?

A simple heuristic argument on the Kac model gives an indication on $j(N)$.
Let us consider a tagged group of $j$ particles and consider,
for any $i = 1,2, \cdots, j$, the set $B_i$ of particles influencing the dynamics of particle $i$ up to the time $t$.
We can assume that the cardinality of $B_i$ is finite to have a correct kinetic behaviour in the limit.
If the propagation of chaos takes place, the groups $B_i$ must be disjoint. On the other hand the probability that two fixed particles 
interact is $O(1/N)$. 
Therefore
the probability that any pair of particles, in the group of $j$, is dynamically correlated is $O(j^2 / N )$
and hence it suffices that  $j/\sqrt{N} \to 0 $  to ensure that the correlations are vanishing. 

The purpose of the present paper is to prove this property for  a class of mean field models including the ones listed above.
One can deal with them simultaneously in terms of an abstract formalism which will be introduced in Section \ref{model}.
Indeed $f_j^N(\calv_j,t)$,  $f_j^N(X_j \calv _j,t)$ and $\r_j^N(t) $ are ruled out by a hierarchy of equations with the same structure
(BBGKY hierarchy). Under suitable hypotheses, the operators occurring in these hierarchies satisfy the same bounds. 
Notice that, physically, these models are very different. 
In particular, the Kac and the soft spheres models are stochastic processes for interacting particles, the quantum mean field
is  time reversible. 

\medskip
Inspired from \cite{PS}, our main tool is a precise notion of decorrelation. Let us present it here,
for the sake of concreteness, in the case of Kac model. 
(For the general definition    in the abstract setting, see Definition \ref{defej} below.) Define
\be
\label{E1}
E_j (t) := \sum_{K \subset J} (-1)^{|K|}  f^N_{J \backslash  K} (t) f(t)^{\otimes K}
\ee
where $J=\{ 1,2, \cdots, j \}$, $K$ is any subset of $J$, $J \backslash K$ 
is the relative complement of $K$ in $J$ and $|K| =$ cardinality of $K $.
$f_A ^N(t)$ stands for the $|A|-$marginal $f_{|A|} ^N(t)$ computed in the configuration $ \{v_i \}_{i\in A} $. 
Similarly, $f(t)^{\otimes K} = f(t)^{\otimes |K|}$ evaluated in $ \{v_i \}_{i\in K}$.
Here $f$ is the solution to the homogeneous Boltzmann equation. 

Eq.\,\eqref{E1} has an inverse formula (proven below), that is
\be
\label{E2}
f^N_j (t) = \sum_{K \subset J}  E_{J \backslash K}(t)  f(t)^{\otimes K}\;,
\ee
where $E_{J \backslash K}$ is defined as $F^N_{J \backslash K}$ is.

$E_j(t)$ measures the tendency of $f^N_j$ to factorize and to converge in the mean field limit. For this reason
$E_j$ will be called the {\it correlation error} (of order $j$) of the mean field system.

These quantities have been already used (under the name ``$v$-functions'') to deal with kinetic limits of stochastic models
\cite{ DP, CDPP91,  BDP,  DOPT,  DOPT1, CP, CPW,  DPTV} and they have been recently investigated in the more singular low density limit of hard spheres \cite{PS}.
In the latter reference, the correlation error of the $N-$particle system is given by
\be
\label{E3}
E'_j (t) := \sum_{K \subset J} (-1)^{|K|}  f^N_{J \backslash K}f^N_1(t)^{\otimes K}\;,
\ee
which is closely related to the cumulant expansion of a probability distribution of particles at a given time.
Note that $E'_j$ 
quantifies the mere deviation of $f^N_j$ from the product of one-particle marginals without any reference to the kinetic equation,
in contrast with \eqref{E1} which measures both factorization and convergence. Unlike in \cite{PS}, in the context of the present paper
\eqref{E1} and \eqref{E3} provide equivalent information, and it is convenient to work directly on the functions
$E_j$ since they satisfy a simple evolution equation.

\bcr
{
Let us make some further comparison with \cite{PS}, where a worst (non optimal) estimate on $E_j$ is obtained. The hard sphere BBGKY hierarchy poses considerably different problems. First of all, in the present paper the analysis is based on the ``correlation equations''. These are driven by a nonlinear mean field problem which is globally well posed, at variance with the inhomogeneous Boltzmann equation. In the setting of \cite{PS}, there is no simple analogue of such correlation equations (notice that Equation \eqref{cij} below fails in this setting, together with the algebraic derivation in the Appendix). As a consequence, in \cite{PS} only the hierarchy for $f^N_j$ is used and a direct expansion to reconstruct and estimate $E_j$ (going through the  definition \eqref{E3} as an intermediate step). Moreover, and most importantly, in \cite{PS} the dynamical correlations are due to collisions which are strong and localized, but rare (`recollisions'). In particular, for hard spheres in the Boltzmann--Grad limit, the interaction  operator (`$T_j/N$' in Equation \,\eqref{eqhiera} below) is not small, and it is replaced by suitable boundary conditions on collision surfaces of diameter $1/\sqrt N$. The smallness of the recollisions is therefore a problem of geometrical nature.
}

The plan of the paper is the following. In Section \ref{modmain} we introduce the model, state our main results 
and write the correlation equations satisfied by $E_j$. The proofs of the results are presented in Section \ref{proofmain},
while the derivation of the correlation equations is given in the Appendix. Similar equations for the evolution of the correlation error have been derived in many of the aforementioned references for stochastic systems.
Finally, Section \ref{applications} collects comments on applications of the results and comparisons with the existing literature.

Let us  conclude this introduction with a remark on the fundamental case of classical particles, i.e.\,the mean field limit of a Newtonian system. 
This case eludes our abstract setting and strategy. Indeed the presence of derivatives makes singular the BBGKY operator, which would
require the introduction of analytic spaces (see \cite{GPP}). 
An efficient approach for the convergence to the Vlasov equation is the direct control of the empirical measures,
and the problem is naturally solved in weak topologies \cite{NeunWick, BraunHepp, Dobru,HaJa}. Concerning the size
of chaos we refer to \cite{GMP16}, where it is shown (avoiding empirical measures) that the $j-$marginal converges to the factorized state with rate 
bounded by $\big(j/N \big)^{1/p}$, for any $p$-Wasserstein distance with
$p \geq2$, while the bound is $j/ \sqrt{N}$ for $p=1$.

\section{Abstract model and main result}\label{modmain}

\subsection{The model}\label{model}

In this section we introduce an abstract setting which extends the trace-class operator formalism customary 
in quantum mechanics. This allows to deal, simultaneously, with classical cases. 
First, we assume the existence of states admitting a family of marginals (Sec.\,\ref{ss:SPS}, \ref{ss:M}),
then we state assumptions on the evolution operators and write down the evolution equations
(Sec.\,\ref{ss:EE}, \ref{ss:H}) and finally we introduce the correlation error (Sec.\,\ref{ss:CE}).

\subsubsection{State of the particle system} \label{ss:SPS}

 Let $\mathfrak{H}$ be a complex separable Hilbert space and let $\mathfrak{H}_n$, $n=1,2, \cdots,N$
 be the $n-$fold tensor power of $\mathfrak{H}$ 
 $$
 \mathfrak{H}_n:=\mathfrak{H}^{\otimes n},  \qquad \mathfrak{H}_1 :=\mathfrak{H}\;.
$$
For any operator $A$ acting on $ \mathfrak{H}_n$, we denote by $A^*$ its adjoint.  As usual $|A|= \sqrt {A^*A}$.
To unburden notations, we drop the $n-$indices here an below, when no confusion arises.

For $1 \leq k,s \leq n$, let $\sigma_{k,s}= \sigma_{s,k}\in  {\cal B} ( \mathfrak{H}_n)$ be the unitary, involutory operator 
defined by
$$
\sigma_{k,s} (e_1 \otimes \cdots \otimes e_k\otimes \cdots \otimes e_s \otimes \cdots e_n)=
e_1 \otimes \cdots \otimes e_s\otimes \cdots \otimes e_k \otimes \cdots e_n\;,\ k\neq s
$$
for any orthonormal basis $\{e_j\}_{j \geq1}$ of $ \mathfrak{H}$.
By convention, we define $\sigma_{k,k}=\mathbb I_{\mathfrak{H}_n}$.
A given $A \in {\cal B} ( \mathfrak{H}_n)$ is said symmetric if
$$
A=\s_{k,s}\,A \,\s^*_{k,s}
$$
for all $k,s = 1,\cdots, n$. 

For any $n$, we postulate the existence of a *-algebra $\mathfrak{A}_n \subset {\cal B}(\mathfrak{H}_n)$,
$\mathfrak{A}_n = \mathfrak{A}^{\otimes n}$, $ \mathfrak{A}_1 =  \mathfrak{A}$
(possibly not containing the identity),  stable under the map $A \to | A |$. 

We assume $\mathfrak{A}_n$ to be equipped with a norm $\norm{\cdot}_1$ defined  by
$$\|A\|_1 := \ell(|A|) \,\,\,\,\,\,\,\,\,\,\,\,\,\,\,\,\,\, A \in \mathfrak{A}_n\;,$$
where
$$
\ell: \mathfrak{A}_n  \to \bC
$$
is a positive linear functional satisfying the following properties:
\begin{itemize}
\item $ |\ell (A)| \leq \ell(|A|)$,\ \ \ $A \in \mathfrak{A}_n$;
\item  for $A \in  \mathfrak{A}_{j} $ and $B \in  \mathfrak{A}_{n-j} $ 
$$
\ell(A \otimes B)=\ell (A) \ell (B);
$$
\noindent as a consequence
$$
\| A \otimes B \|_1 = \ell (| A \otimes B |) = \ell (| A |\otimes |B |) = \| A \| _1 \| B\|_1;
$$
\item for any $A \in  \mathfrak{A}_{n} $,  and any $\sigma_{k,s}$,
$$
\| \sigma_{k,s} A  \sigma_{k,s}^* \|_1= \| A \|_1. 
$$
\end{itemize}

We consider the completion of the above algebras and keep for it the same notation.
Thus $\AA_n$  with norm $\|\cdot\|_1$ is a Banach space
and we extend $\ell$ as a continuous functional by completion.

A state of the $N-$particle system is, by definition, an element $F^N  \in \mathfrak{A}_N$,
positive (as operator in  ${\cal B} ( \mathfrak{H}_N)$), symmetric and such that $\| F^N \|_1 = 1$.

\subsubsection {Marginals} \label{ss:M}
Let
$\ell_n: \mathfrak{A}_n \to \mathfrak{A}_{n-1}$ be positive linear maps such that:
\begin{itemize}
\item $\ell_n(A_1 \otimes \cdots \otimes A_n) = \ell(A_n) (A_1 \otimes \cdots \otimes A_{n-1})\;,\ \ \ \ A_i \in \AA\;.$
\end{itemize}
Note that $\ell_n$ is symmetry preserving.

Moreover we assume that
\begin{itemize}
\item $ |\ell_n (A)| \leq \ell_n(|A|)$,\ \ \ $A \in \mathfrak{A}_n$.
\end{itemize}
\vskip 0.3cm

For $1 \leq j < n $, we indicate by $\ell_{j,n}:\AA_n \to \AA_j$ the transformation
$$\ell_{j,n} =\ell_{j+1}\cdots \ell_{n}\;.$$
This map is a contraction and  preserves the norm of positive elements:
\be
\label{contr}
\| \ell_{j,n} (A) \|_1 = \ell(|\ell_{j,n} (A)|) \leq \ell(\ell_{j,n}(|A|)) =   \| A \|_1
\ee
for $A \in \AA_n$ and the equality holds for $A$ positive.

The $j-$particle marginal of the $N$-particle state $F^N$ is given by
$$
F^N_j :=\ell_{j,N} \left(F^N\right) \in \mathfrak{A}_j .
$$
By construction, $F^N_j$ is a $j$-particle state.


\subsubsection {Evolution equations} \label{ss:EE}

The evolution of a state $\FF^N \to \FF^N(t)$ in $ \mathfrak{A}_N$ is supposed to be given by the $N-$particle dynamics associated to a two-body interaction:
\be\label{Neq}
\frac{d}{dt}\FF^N=(K_0^N+V^N)\FF^N,
\ee
where 
 \be\label{defV}
K_0^N=\sum_{i=1}^N \mathbb I_{\mathfrak{A}_{i-1}}\otimes K_0\otimes\mathbb I_{ \mathfrak{A}_{N-i}}
\ee 
and
\be\label{defV1}
V^N=\frac{1}{N}\sum_{1\leq i< j\leq N}V_{i,j},
\ee
with
\be\label{defV}
V_{i,j}(A):=\sigma^*_{1,i}\sigma^*_{2,j}V
\otimes  \mathbb I_{ \mathfrak{A}_{N-2} }(\sigma_{1,i}\sigma_{2,j}A
\sigma^*_{2,j}\sigma^*_{1,i})\sigma_{2,j}\sigma_{1,i}, \ A\in \mathfrak{A}_{N}
\ee
for a (possibly unbounded) linear operator $K_0$ on $\AA$ and a symmetry preserving, two-body potential $V\in
{\cal B}(\AA_{2})$. 

Formula \eqref{defV} 
expresses the following simple fact. If
$A=A_1 \otimes \cdots \otimes  A_N $ and
$$
V(A_i \otimes A_j) = \sum_{r,s} \alpha_{r,s} B_r \otimes C_s\;,
$$
then
$$
V_{i,j} (A)= \sum_{r,s} \alpha_{r,s} A_1 \otimes \cdots A_{i-1} \otimes B_r\otimes  \cdots
\otimes  C_s \otimes A_{j+1} \cdots \otimes A_N \;.
$$
In other words, the action of $V_{i,j}$ is the same as  $V_{1,2}=V$  on the  slots $i$ and $j$ 
(see the table below for concrete examples).

We assume that both $K_0$ and $K_0^N + V^N $ generate a strongly continuous, positivity preserving, isometric semigroup (with respect to the norm $\|\cdot\|_1$)
and there exists a unique mild solution to \eqref{Neq} with initial datum $F^N(0) \in \AA_N$. 
Symmetry is automatically preserved by the symmetry of $ K^N_0$ and $ V^N$.

Finally, for any $F\in \mathfrak{A}$, $F^N \in  \mathfrak{A}_N $ and $i,r >j$, $i\neq r$  we assume
\be\label{trackf}
\ell\left(K_0\left(F\right)\right)=0\mbox{ \ \ \ and\ \ \  }
\ell_{j,N} \left(V_{i,r} \left(F^N\right)\right)=0\;.
\ee
These last properties are necessary to deduce the forthcoming hierarchy.

The following table summarizes the applications of the above abstract model to the three models presented in the introduction, namely the {\bf K}(ac), {\bf S}(oft spheres) and {\bf Q}(uantum) models. \bcb The precise settings and statements will be given below, in Section \ref{applications}.\ec  For these models, we have
that $\mathfrak{H}$ is $L^2(\bR^3,dv), L^2(\bR^{6},dxdv)$ and $L^2(\bR^d,dx)$ respectively, while
$\mathfrak{A}$ is given by $L^1 \cap L^\infty (\bR^{3}), L^1 \cap L^\infty (\bR^{6})$ (as multiplication operators\footnote{The
restriction to $L^\infty$ is merely due to the abstract formulation. This assumption 
can be removed by density in the main theorem \ref{main} below.})
and ${\cal L}^1(L^2(\bR^d))$ (the space of the trace-class operators).

\bigskip

\begin{center}
   \begin{tabular}{ | r | c | c |c|  }
     \hline
     &Kac&Soft spheres&Quantum mean field\\
     \hline
  $\ell $&  $f \to \int dv f $&  $f \to \int \int f dx dv $& $A \to \tr A$\\
  \hline
  $\|\cdot\|_1$&  $\|\cdot\|_{L^1}$&  $\|\cdot\|_{L^1} $& $A \to \tr |A|$
  \\
  \hline
  $K_0 $&   $0$&  $-v \cdot \nabla_x $&  $ \frac1{i\hbar}[ -\frac{\hbar^2}{2}\D, \cdot  ] $\\
  \hline
  $V_{i,j}(F^N)$ &  \begin{tabular}{c}
  $\int d\omega$\\$ B(\omega; v_i-v_j)$\\$\{f^N(\calv_N^{i,j})$\\$-f^N(\calv_N)\} $\end{tabular}
& \begin{tabular}{c}$h\left(|x_i-x_j|\right) $\\$B\left(\frac{x_i-x_j}{|x_i-x_j|}; v_i-v_j\right)
$\\$\{f^N(X_N ,\calv_N^{i,j})$\\$-f^N(X_N,\calv_N)\} $\end{tabular}
& $ \frac1{i\hbar}[ V(x_i-x_j), \rho^N ]$\\
\hline
     \end{tabular}
     \end{center}
 
 \bigskip  
 
    
\noindent
In both cases $\bS$ and $\bQ$, $ K_0$ is only densely defined. As we shall see, we will use only the groups generated by $ K_0$, i.e.\,respectively
$$
e^{K_0t}f (x,v)=f(x-vt,v)\;, \quad f \in L^{1} \cap L^\infty
$$
for the case $\bS$, and
$$
e^{K_0t}A = U_0(-t)AU_0(t)\;, \qquad A\in {\cal L}^1
$$ 
where $U_0(t)=  e^{ -i \frac{\hbar}{2}\Delta_x t}$ as unitary operator on $L^2$ for the case $\bQ$.
The operators $V_{i,j}:  \mathfrak{A}_N \to  \mathfrak{A}_N $  have to be understood as bounded operators from
$L^{\infty} \cap L^1 \to L^{\infty} \cap L^1$
 for both cases $\bK$ and $\bS$.

\subsubsection{Hierarchies}  \label{ss:H}

Applying subsequently $\ell_{N}, \ell_{N-1},\cdots$ to \eqref{Neq} and using the symmetry and Eq.\,\eqref{trackf}, we get the BBGKY hierarchy of equations
\be\label{eqhiera}
\pa_t \FF^N_j = \left(K_0^j+\frac {T_j}{N}\right) \left(\FF^N_j\right)+ \a(j,N) C_{j+1}\left( \FF^N_{j+1}\right)
\ee
for $j=1,\cdots,N$, where 
\be
\a(j,N) = \frac {(N-j)} N\;,
\label{notalp}
\ee
$K_0^j, T_j$  are operators on $\mathfrak{A}_{j}$,
 \be
K_0^j:=\sum_{i=1}^j \mathbb I_{ \mathfrak{A}_{i-1}}\otimes K_0\otimes\mathbb I_{\mathfrak{A}_{j-i}}\;,
\ee
\be
T_j:=\sum_{1\leq i < r \leq j} T_{i,r}
\ee
with\footnote{According to our definition, it should be specified that $T_{i,r}:\ \mathfrak{A}_{j}\to \mathfrak{A}_{ j} $ depends explicitly on $j$. We avoid to introduce a further notation, this fact being  clear from the context.}
\be
 T_{i,r}(A)=
\sigma^*_{1,i}\sigma^*_{2,r}V
\otimes  \mathbb I_{ \mathfrak{A}_{j-2} }(\sigma_{1,i}\sigma_{2,r}A
\sigma^*_{2,r}\sigma^*_{1,i})\sigma_{2,r}\sigma_{1,i}\;, \,\,\,\,\, A\in \mathfrak{A}_{j}
\ee
and
\be
\label{C}
C_{j+1}\left(A\right):=\ell_{j+1}\left(\sum_{i\leq j}T_{i,j+1}\left(A\right)\right)=\sum_{i=1} ^jC_{i,j+1}\left(A\right)\;, \,\,\,\,\, A\in \mathfrak{A}_{j+1}
\ee
with
\be
C_{i,j+1}\left(A\right):=\ell_{j+1}\left( T_{i,j+1}\left(A\right) \right)\;,
\label{eq:Ci}
\ee
\be
C_{i,j+1}:  \mathfrak{A}_{j+1}\to \mathfrak{A}_{ j}.
\ee
We have
\be\label{normtc}
\hspace{0,3cm}\norm{T_j}\leq \frac {j(j-1)}2 \| T_{i,r}\| \leq  \frac {j(j-1)}2 \| V\|\;,
\hspace{0,8cm} \norm{C_{j+1}}\leq j \|C_{i,j+1}\|\leq j \|V\|\;
\ee
(meant for $\norm{T_j}_{\mathfrak{A}_j\to\mathfrak{A}_j},\ \norm{T_{i,r}}_{\mathfrak{A}_j\to\mathfrak{A}_j},\ \norm{C_{j+1}}_{\mathfrak{A}_{j+1}\to\mathfrak{A}_{j}},\ \norm{C_{i,j+1}}_{\mathfrak{A}_{j+1}\to\mathfrak{A}_{j}},\ \norm{V}_{\mathfrak{A}_2\to\mathfrak{A}_2}$).

Associated to $V$, we introduce the nonlinear mapping $Q: \mathfrak{A}^{\otimes 2}  \to \mathfrak{A}$ by the formula
\be\label{nlc}
 Q(\FF,G):=\ell_2(V(\FF\otimes G))
\ee
and the nonlinear mean field equation on $\mathfrak{A}$
\be\label{mfe}
\pa_t\FF= K_0 \left(\FF\right)+Q(\FF,\FF),\ \ \ \ F(0)\geq 0,\ \ \ \ \norm{F(0)}_1=1\;.
\ee
Eq.\,\eqref{mfe} is the Boltzmann, Povzner or Hartree equation according to the specifications
established in the table above. 
By assumption, it possesses a (global) unique mild solution, preserving $\|\cdot\|_1$ and
positivity.

Observe that, after definitions \eqref{nlc} and \eqref{eq:Ci},
\be
\label{cij}
C_{i,j+1}\left(\FF^{\otimes(j+1)}\right)=\FF^{\otimes(i-1)} \otimes Q(\FF,\FF) \otimes \FF^{\otimes(j-i)}
\ee
and, if $F$ is a solution of \eqref{mfe}, then $\{\FF^{\otimes j}\}_{  j \geq 1}$ solves
\be\pa_t\FF^{\otimes j}= K_0^j \left(\FF^{\otimes j}\right)+C_{j+1}\left(\FF^{\otimes(j+1)}\right)\;.
\ee
In other words, $\{\FF^{\otimes j}\}_{  \ j \geq 1} $ is  a solution of the formal limit of the hierarchy \eqref{eqhiera} as $N\to\infty$ .
%


\subsubsection{Correlation error.} \label{ss:CE}

Let us now fix $j$ and $K\subset\{\un{j}\} = J$. Writing $K=\{i_s; s=\un{k}\}$ with $i_1 < i_2 < \cdots <i_k$,
$|K| = k$, we define the unitary operator
\be
\sigma_{K,j} =
\sigma_{i_1,1} \sigma_{i_2,2}\cdots \sigma_{i_k,k}\;,
\ee
$\sigma_{\emptyset,j} = \mathbb I_{\mathfrak{H}_j}$. 
We introduce a mapping 
on $A_{k} \in \mathfrak{A}_{k}$   
into $ \cal {B} ( \mathfrak{H}_{j})$
\be
A_{k} \longrightarrow A^J_{K}\;,
\label{o}
\ee
by 
\be
\label{AJ}
A^J_{K}= \s_{K,j} (A_{k} \otimes \mathbb I_{\mathfrak{H}_{j-k}}) \s_{K,j}^*\;.
\ee

For instance in the case of marginals, dropping from now on the explicit dependence on $N$ ($F^N_j=F_j$), we get
\be
F^J_{J \backslash K}= \s_{J \backslash K,j} (F_{j-k} \otimes \mathbb I_{\mathfrak{H}_{k}}) \s_{J \backslash K,j}^*\;.
\ee

Moreover, for $F \in \AA$ we write
\be
\label{ODt}
F^{\otimes K,J} = \s_{J \backslash K,j} (\mathbb I_{\mathfrak{H}_{j-k}} \otimes F^{\otimes k} ) \s_{J \backslash K,j}^* \in \cal {B} ( \mathfrak{H}_{j}) \;.
\ee
A more explicit equivalent definition is
$$
F^{\otimes K,J}=A_1 \otimes A_2\otimes \cdots \otimes A_j\;,
$$
$$
A_i=F \,\,\,\text {if} \,\,\,\,\,\, i \in K  \,\,\,\, \text {and}\,\,\,\,\,\, A_i=\mathbb I_{\mathfrak{H}} \,\,\,\,\, \text {otherwise.}
$$
Note that 
\be
F^{\otimes K,J}F^{\otimes R,J}=F^{\otimes R,J}F^{\otimes K,J} =F^{\otimes (K\cup R),J}
\label{eq:EI}
\ee
if $R\cap K =\emptyset$ both in $J$. Note also that $\FF^{\otimes K,J }$ and $A^J_{J \backslash K}$
commute.

%

Even if $\mathbb I_{\mathfrak{H}_{k}}\notin\mathfrak{A}_{k}$ (and $A^J_{J \backslash K} \notin \mathfrak{A}_{j}$,
$F^{\otimes K,J} \notin \mathfrak{A}_{j}$), a product of the form $\FF^{\otimes K,J } A^J_{J \backslash K}$ lies always in $\AA_j$ and
\be
\|\FF^{\otimes K,J } A^J_{J \backslash K}\|_1 = \|F\|^k_1\|A_{j-k}\|_1\;.
\label{eq:NppE}
\ee

Fixing the convention
\be
F^{\otimes \emptyset,J} =  A^J_{\emptyset} = \mathbb I_{\mathfrak{H}_j}
\label{eq:con0J}
\ee
we introduce, in the following definition, a family of symmetric elements in $\mathfrak{A}_j$ characterizing the state of the $N-$particle system.

\newcommand{\Pii}[2]{\underset{#1}{\overset{#2}{\Pi}}}
\begin{Def}\label{defej}
For any 
$j=\unN$, setting $J=\{\un{j}\}$ and $k = |K|$,
we define the  
``correlation error" of order $j$ by
\be\label{deferror}
E_j=\sum_{K \subset J} (-1)^{k} \ 
\FF^{\otimes K,J} \FF^J_{J \backslash K}\;.
\ee
\end{Def}

\noindent
Eq.\,\eqref{deferror} also reads
\be\label{deferrorabs}
E_j=\sum_{K \subset J} (-1)^{k} \ 
\s_{J \backslash K,j}(\FF_{j- k} \otimes \FF^{\otimes k})\s_{J \backslash K,j}^*\;.
\ee
Note that, by \eqref{eq:con0J}, the terms $K=\emptyset$ and $K = J$ have to be interpreted 
as $F^J_J = F_j$ and $(-1)^j F^{\otimes J,J} = (-1)^j F^{\otimes j}$ respectively.

Formula \eqref{deferror} is inverted by
\be\label{invdeferror}
F_j=\sum_{K \subset J} F^{\otimes K,J} E^J_{J \backslash K} 
\ee
where, using \eqref{o}, we write
\bea
E^J_{J \backslash K}&=&\s_{J \backslash K,j} (E_{j-k} \otimes \mathbb I_{\mathfrak{H}_{k}}) \s_{J \backslash K,j}^*
\label{eq:dakaK}\\
&=& \sum_{R \subset J\backslash K} (-1)^{|R|}  \FF^{ \otimes R,J}  \FF^J_{J \backslash (K\cup R)}
\label{eq:dakaKb}
\eea
and $E^J_{\emptyset} = \mathbb I_{\mathfrak{H}_j}$.
To prove \eqref{eq:dakaKb}, denoting $I = \{1,\cdots,j-k\}$ and $J \backslash K=\{i_s; s=\un{j-k}\}$ with $i_s$ increasing, 
we observe that \eqref{deferrorabs} together with the change of variables induced by $\s_{J \backslash K,j}=
\sigma_{i_1,1} \cdots \sigma_{i_{j-k},j-k}$,
$$I \supset R'=\{\ell_1,\cdots,\ell_{|R|}\} \to R=\{i_{\ell_1},\cdots,i_{\ell_{|R|}}\} \subset J \backslash K\;,$$
imply
\be
 E^J_{J \backslash K} = \sum_{R \subset J\backslash K} (-1)^{|R|} \s_{J \backslash K,j} (
\s_{I \backslash R',j-k}(\FF_{j- k-|R|} \otimes \FF^{\otimes |R|})\s_{I \backslash R',j-k}^*
 \otimes \mathbb I_{\mathfrak{H}_{k}}) 
 \s_{J \backslash K,j}^*\;. \nn
 \ee
 On the other hand
 \bea
&& \s_{J \backslash K,j} (
\s_{I \backslash R',j-k}(\FF_{j- k-|R|} \otimes \FF^{\otimes |R|})\s_{I \backslash R',j-k}^*
 \otimes \mathbb I_{\mathfrak{H}_{k}}) 
 \s_{J \backslash K,j}^*\nn\\
&&= \s_{J \backslash K,j} 
(\s_{I \backslash R',j-k}\otimes I_{\mathfrak{H}_{k}})( \FF_{j- k-|R|} \otimes\FF^{\otimes |R|}\otimes\mathbb I_{\mathfrak{H}_{k}})
(\s^*_{I \backslash R',j-k}\otimes I_{\mathfrak{H}_{k}})
 \s_{J \backslash K,j}^*\nn\\
 &&=  \s_{J \backslash K,j} 
(\s_{I \backslash R',j-k}\otimes I_{\mathfrak{H}_{k}})
(\mathbb I_{\mathfrak{H}_{j-k-|R|}} \otimes F^{\otimes |R|}\otimes  \mathbb I_{\mathfrak{H}_{k}})
(\s^*_{I \backslash R',j-k}\otimes I_{\mathfrak{H}_{k}})
\s^*_{J \backslash K,j} 
\nn\\
&&\ \ \ \cdot  \s_{J \backslash K,j} 
(\s_{I \backslash R',j-k}\otimes I_{\mathfrak{H}_{k}}) (F_{j-k-|R|} \otimes \mathbb I_{\mathfrak{H}_{k+|R|}})
(\s^*_{I \backslash R',j-k}\otimes I_{\mathfrak{H}_{k}})
 \s_{J \backslash K,j}^*\nn\\
 &&=
 \s_{J \backslash R,j} (\mathbb I_{\mathfrak{H}_{j-|R|}} \otimes F^{\otimes |R|} ) \s_{J \backslash R,j}^*\nn\\
 &&\ \ \ \cdot \s_{J \backslash (K\cup R),j} (F_{j-k-|R|} \otimes \mathbb I_{\mathfrak{H}_{k+|R|}}) \s_{J \backslash (K\cup R),j}^*\;.
 \nn
 \eea
Formula \eqref{invdeferror} follows then from \eqref{eq:dakaKb} and \eqref{eq:EI}:
\bea
\sum_{ K \subset J} F^{ \otimes K,J} E^J_{J \backslash K} & =& \sum_{K \subset J} \sum_{R \subset J \backslash K} (-1)^{|R|} F^{\otimes K,J} F^{\otimes R,J} F^J_{J \backslash (K \cup R)} \nn \\
& =& \sum_{K \subset J} \sum_{R \subset J \backslash K} (-1)^{|R|} F^{\otimes (K \cup R),J} F^J_{J \backslash (K \cup R)} \nn \\
& = &\sum_{T \subset J} \sum_{R \subset T} (-1)^{|R|} F^{\otimes T,J} F^J_{J \backslash T}\nn\\
& = &\sum_{T \subset J} \delta_{T,\emptyset} F^{\otimes T,J} F^J_{J \backslash T} \nn\\
& =& F^J_J \equiv F_j\;. \nn 
\eea

\medskip

\subsection{Main result}\label{mainresult}

Let $\FF^N(t)$ be the  time-evolved state of the $N-$particle system, solution of \eqref{Neq} with initial datum $\FF^N(0)\in \AA_N$. Let $\FF(t)$ be the solution of the kinetic equation \eqref{mfe} with initial datum $\FF\in \AA$ and $E_j(t), j=\unN,$ the correlation errors associated to the marginals of $\FF^N(t)$, as given by Definition \ref{defej}.

In the sequel, we will denote  for any operator $H:\mathfrak{A}_{j}\to\mathfrak{A}_{j'},\ j,j'=1,\dots, N$, 
\be\label{abuse}
\|H\| = \|H\|_{\mathfrak{A}_{j}\to \mathfrak{A}_{j'}}.
\ee
\begin{Thm}\label{main}
 Let us suppose  
that, for all $j=1,\dots, N $ and for some $C_0\geq 1$,
\be
\label{hyp}
\norm{E_j(0)}_1\leq C_0^j\left(\frac{j}{\sqrt N}\right)^{j}.
\ee
Then,
 for all  $t>0$ and all $j=~1,\dots, N$, 
 one has 
\be
\label{eqmain}
\| E_j(t) \|_1 \leq \left( C_2 e^{C_1 t\norm{V}} \right)^j \left(\frac{j}{\sqrt N}\right)^j 
\ee
where $C_1 \geq 0,\ C_2\geq 1 $ are 
given by  formula \eqref{c1c2}.

Suppose in addition that $\norm{E_1(0)}_1\leq \frac {B_0}N$ for some $B_0>0$. Then for all $t > 0$
\be\nonumber
\norm{E_1(t)}_1\leq \frac 1N \big( B_2 e^{B_1 t \|V\| } \big)\;
\ee
for   $B_1>0,\ B_2\geq 1 $  
given by formula \eqref{b1b2}.
\end{Thm}

\begin{Cor}\label{cormain}
Suppose that \eqref{hyp} holds and that $\norm{E_1(0)}_1\leq \frac {B_0}N$ for some $B_0>0$. Then for all  $t>0$ and all $j=~1,\dots, N$,  the marginals satisfy
 \be\label{eqcormain}\nn
\| \FF_j (t) - \FF(t)^{\otimes j}   \|_1 \leq 
 D_2 e^{D_1t\norm{V}}\frac {j^2 }N
\ee
where $D_2=\sup\{B_2,{8(eC_2)^2}\},\  D_1=\sup\{B_1,2C_1\}$.
\end{Cor}

\bcrr{
\begin{proof} 
We have, according to \eqref{invdeferror} and Theorem \ref{main}, 
   
\bea
&\| F_j (t) - F(t)^{\otimes j}   \|_1 \nn\\
&  \leq   \sum\limits_{k=1}^j \binom {j}{k} \| E_k(t) \|_1
\leq \| E_1(t) \|_1 +   \sum\limits_{k=2}^j \binom {j}{k} (\frac {k^2  C_2^2e^{2C_1 t\norm{V}}}N )^{k/2} \nn \\
& =  \| E_1(t) \|_1 + \sum\limits_{k=2}^j j(j-1) \cdots (j-k+1) \frac {k^k}{k!} (\frac { C_2 e^{C_1t\norm{V}}} {\sqrt {N}})^k  
 \nn \\
& \leq  \| E_1(t)\|_1+ \frac1{\sqrt{2\pi}}\sum\limits_{k=2}^j (\frac {ej  C_2 e^{C_1t\norm{V}}}{ \sqrt {N}})^k 
\leq 
\| E_1(t) \|_1 +\frac1{\sqrt{2\pi}} \frac { \big(\frac {ej  C_2e^{C_1t\norm{V}}}{ \sqrt {N}}\big)^2}{1-
\frac {jeC_2 e^{C_1\norm{V}t}}{ \sqrt {N}}}\nn
\eea
since $\frac{k^k}{k!}\leq \frac{e^k}{\sqrt{2\pi k}}$.
Then the result follows if 
\be
\label{Nlarge}
N \geq 4(jeC_2)^2 e^{2C_1\norm{V} t}.
\ee
 On the other hand if \eqref{Nlarge} is violated 
 $$
 8 \frac { (ejC_2)^2 e^{2C_1t\norm{V}}}N \geq \frac 84=2\geq \| F_j (t) - F(t)^{\otimes j}   \|_1
 $$
 since $F_j (t)$ and $F(t)^{\otimes j}  $ remain normalized for all $t$, and this concludes the proof.
 }
\end{proof}

%
%
%
%
%

\begin{Rmk}
The bound \eqref{eqmain} is trivial when $j \geq \sqrt{N}$, by virtue of the obvious inequality $\| E_j (t) \|_1 \leq 2^j$.
Therefore we will consider in the sequel only the case $j < \sqrt{N}$.
Furthermore note that there is no need in Corollary \ref{cormain} for the initial condition $\FF^N(0)$ to be factorized. 
\end{Rmk}

\begin{Rmk}
\bcr
{
It may be worth  discussing the meaning of the hypothesis that  $\|V\|$ is bounded in the three concrete models described in the introduction. For the Kac model, as well as for soft spheres, this boundedness of $\|V\|$ means that the cross--section for the associated Boltzmann equation must be bounded, as required e.g. in \cite{GM97}. This condition is often referred as ``pseudo-Maxwellian cross--section''. From a physical point of view, particles interact via a specific inverse power law potential, and an angular cutoff is also applied. However, beyond this case, there are many physically interesting situations fulfilling the boundedness condition. 
An example is the quantum Boltzmann equation which has a similar form as the classical Boltzmann equation. In this case the cross-section is bounded, provided the interaction potential has suitably decaying Fourier transform \cite{BCEP04}.
}

\bcr
{
Unfortunately we do not handle here more general cross--sections diverging with the relative velocity as, for instance, the hard--sphere model. The hierarchical approach does not seem to work conveniently in this case. For example in  \cite{ACI91}, in order to obtain a uniqueness result on the solutions of the hard sphere hierarchy, it is made use of the equivalence with the notion of statistical solutions; see also \cite{MM}.}

\bcr
{
For the quantum mean--field regime the boundedness condition is a simple consequence of the requirement that the interaction potential is bounded.  
}
\end{Rmk}

 Theorem \ref{main} will be proven in Section \ref{proofmain}.

\subsection{The correlation equations}\label{equa}
In this section we write the equations satisfied by the errors $E_j$  introduced in Definition \ref{defej}.



We make use of the notation \eqref{notalp} and $J^{i} = J \backslash \{i\}$, $J^{i,s} = J \backslash \{i,s\}$.
We introduce four (time-dependent) operators $D_j,D_j^1,D_j^{-1}$ and $D_j^{-2}$, $j=1,\cdots, N$, by the following formulas.

\bea
D_j : \ \AA_j&\to&\AA_j\nn\\
A_j &\mapsto& \a(j,N)  \suml_{ i\in J}  C_{i,j+1}  \left( {\FF} ^{\otimes \{i\},J\cup \{j+1\}} A_{J^{i} \cup \{j+1\}}^{J\cup \{j+1\}} +F^{\otimes \{j+1\},J\cup \{j+1\}} A_J^{J\cup \{j+1\}} \right) \nn \\
&& -\frac 1N  \sum_{\substack{i,s \in J \\ i \neq s}}
  C_{i,j+1} \left(F^{\otimes \{ s\},J\cup \{j+1\}}  A^{J\cup \{j+1\}}_{J^s \cup \{j+1\}} \right) \nn\\
  \label{eqdefDj}
\eea
where we used \eqref{o} and \eqref{ODt} (and the convention $ \suml_{i \neq s \in J}=0$ for $J = \{1\}$). 
Here $F = \FF(t)$ is the solution of \eqref{mfe}. The meaning of the above operator is transparent. Given $A_j$ we can construct via formula \eqref{AJ} the operators $A^{J \cup \{ j+1\}}_S$ with $|S|=j$. The right hand side of the above expression yields an operator in  $\AA_j$. Similar arguments apply in the following formulas.

From now on, to unburden the notation we will drop the upper indices of set, i.e.\,$D_j,\ \AA_{j}\to\AA_j$, is written as
\bea
 A_j &\mapsto& \a(j,N)  \suml_{ i\in J}  C_{i,j+1}  \left( {\FF} ^{\otimes \{i\}} A_{J^{i} \cup \{j+1\}} +F^{\otimes \{j+1\}} A_J  \right) \nn\\
&& -\frac 1N  \sum_{\substack{i,s \in J \\ i \neq s}}
   C_{i,j+1} \left(F^{\otimes \{ s\}}  A_{J^s \cup \{j+1\}} \right) \;. 
   \label{eqdefDj'}
\eea
Analogously, we define:
\bea
D^1_j: \AA_{j+1}&\to&\AA_j,\ \ \ j=1,\dots,N-1,\nn\\
A_{j+1}&\mapsto&
\a(j,N)  C_{j+1} 
\left(A_{j+1}\right)\;,
   \label{eqdefDj1}
\eea
\bea
D^{-1}_j: \AA_{j-1}&\to&\AA_j,\ \ \ j=2,\dots,N,\nn\\
A_{j-1} &\mapsto&
-\frac jN  \suml_{i\in J} Q(\FF,\FF)^{\otimes\{i\}}A_{J^{i}} + \frac 1{N} \suml_{\substack{i,s \in J \\ i \neq s}}  T_{i,s} \left({\FF}^{\otimes\{i\}}A_{J^{i}}\right)\;+ \nn\\ 
&&-\frac 1N  \suml_{\substack{i,s \in J \\ i \neq s}}  C_{i,j+1} \left ( {\FF}^{\otimes\{i,s\}} A_{ J^{i,s} \cup \{j+1\}} +   {\FF}^{\otimes\{s,j+1\}}A_{J^s} \right )\;,
\label{eqdefDj-}
\eea
\bea
D^{-2}_j: \ \AA_{j-2}&\to&\AA_j,\ \ \ j=3,\dots,N,\nn\\
A_{j-2}  &\mapsto&
\frac 1N  \suml_{\substack{i,s \in J \\ i \neq s}} \left ( \frac{1}{2}T_{i,s} 
\left({\FF}^{\otimes\{i,s\}}A_{J^{i,s}}\right) - Q(F,F)^{\otimes\{i\}} 
F^{\otimes \{s\}}A_{J^{i,s}} \right )\;.\label{defd=}
\eea

Let us consider the equation:
\bcrr{
\bea\label{eqhieraerror}
\pa_t E_j&=& \left(K_0^j+\frac{T_j}{N}\right)\left(E_j\right) +
D_j
\left(E_{{j}}\right) \nn \\
&&+ 
D_j^1
\left(E_{j+1}\right) + 
D_j^{-1}
\left(E_{j-1}\right) +
D_j^{-2}
\left(E_{j-2}\right)
\eea
where, by convention, 
\be
\label{E0}
\left\{
\begin{array}{l}
 D_N^1:=D_1^{-2}:=0\\
 D_1^{-1} \left(E_0\right):=-\frac 1N Q(F,F)\;,
\\
D_2^{-2}\left(E_0\right):= \frac 1 N  \big(T_{1,2} F \otimes F -  Q(F,F)^{\otimes\{1\}} 
F^{\otimes \{2\}} -  Q(F,F)^{\otimes\{2\}} 
F^{\otimes \{1\}}\big) \;.
\end{array}
\right.
\ee
}

Note that the first line contains  operators which do not change the particle number. 
$D^1_j$ is an operator increasing 
by one the number of particles considered. $D^{-1}_j$ and $D^{-2}_j$ are operators decreasing the number of particles by one and two respectively.

Eq.\,\eqref {eqhieraerror} is inhomogeneous so that it has nontrivial solutions even for initial data $E_j(0)=0, \, j>0 $ (namely when the initial state is chaotic).

Notice that, by \eqref{normtc}, \eqref{eq:NppE} and the normalization of $F$,
\be
\label{eq:estDj}
\| D_j\| \leq \suml_{ i\in J} 2 \|C_{i,j+1} \| +  \frac 1N  \sum_{\substack{i,s \in J \\ i \neq s}}
\|C_{i,j+1}\| \leq 2 j \|V\| + \frac{j^2}{N}\|V\| \leq 3 j \|V\|\;.
\ee
Similarly,
\be
\label{eq:estD'j}
\| D^1_j\| \leq j\|V\|\;, \quad \| D_j^{-1} \| \leq 4 \frac {j^2} N\|V\|\;,  
\quad \| D_j^{-2} \| \leq \frac{3}{2}\frac {j^2} N\|V\|\;.
\ee

The following result will be proven 
in the Appendix.
\begin{Prop}\label{maineq}
Let $\FF$ satisfy  \eqref{mfe}.
Then the hierarchy of equations \eqref{eqhiera}
 is equivalent to  the hierarchy
\eqref{eqhieraerror}
 in the sense that
\be\nonumber
\{\FF_j\}_{j=1,\cdots, N}\mbox{ solves }\eqref{eqhiera}
\Longrightarrow
\{E_k\}_{k=1,\cdots, N} \mbox{ solves } \eqref{eqhieraerror}\;.
\ee
\end{Prop}

\section{Proof of Theorem \ref{main}}\label{proofmain}
We start by focusing on the evolution generated by the operator $$ K_0^j +\frac{T_j}{N} + D_j$$
preserving the particle number. We construct the (two-parameters) semigroup $U_j(t,s)$ for $s \leq t$, satisfying
\bea\label{u0}
\pa_t U_j(t,s)&=&\left(K_0^j+\frac {T_j}{N}+D_j(t)\right) U_j(t,s),\quad j=1,\dots,N\;,\nn \\
 U_j(s,s)&=& \mathbb{I}_{\mathfrak{A}_j} \;. 
\eea
Recall that 
\be\label{unitu}
\norm{e^{\left(K_0^j+\frac {T_j}{N}\right)t}}=1
\ee
by assumption.
%

Moreover, we assume preliminarily that
\be
\label{normV}
\|V\| = \frac 14.
\ee
From \eqref{normtc} and \eqref{eq:estDj} we get then 
$
\| T_j \|  \leq j^2/8
$
 and 
\be\label{normdj}
\| D_j \| \leq j.
\ee
We deduce, by Gronwall Lemma,
\be\label{gronwald}
\|U_j(t,s) \| \leq e^{j(t-s)}\;.
\ee
We turn next to the contribution of the operators changing the particle number and notice that
(cf.\,\eqref{eq:estD'j})
\be\label{normd}
\| D^1_j\| \leq j, \quad \| D_j^{-1}\| \leq  \frac {j^2} N,  \quad \| D_j^{-2}\|
 \leq \frac {j^2} N\;.
\ee

In order to estimate $E_j$, we express the solution of Eq.\,\eqref{eqhieraerror} in terms of the Dyson series
\bea\label{dysexp}
E_j(t)&= & \suml_{n \geq 0} \suml_{k_1 \cdots k_n} \int_s^t dt_1 \int_s^{t_1} dt_2 \cdots \int_s^{t_{n-1} }dt_n \nn \\
&&U_{j_1}(t,t_1)D_{j_1}^{k_1}(t_1)U_{j_2}(t_1,t_2) \dots D_{j_n}^{k_n} (t_n) U_{j_0}(t_n,s) \left(E_{j_0} (s)\right)\;,
\eea
where:
\begin{itemize}
\item the term $n=0$ is $U_j(t,s)\left(E_j(s)\right)$;
\item $\suml_{k_1 \cdots k_n}$  denotes the sum over all possible choices for the sequence $k_1\dots k_n$ with $k_i \in \{1,-1,-2\} $;
\item $ j_1=j, j_2=j+k_1, \cdots, j_{m+1}=j_m+k_{m}$ and 
$j_0$ is the number of particles at time $s$, namely $j_0=j+ \suml_{\ell=1}^n k_\ell=j_n+k_n $;
\item we use the convention expressed by \eqref{E0}. 

\item  the Dyson series \eqref{dysexp} follows by iterating \eqref {eqhieraerror}  in integral form via the Duhamel formula, and since the terms
$D_N^{1}(E_N)=0$ and $D^{-1}_1\left(E_0\right),$ $D_1^{-2}\left(E_{-1}\right),$ 
$D^{-2}_2\left(E_0\right)$ are explicitly known (see \eqref {E0}), the iteration stops when these terms appear; 
namely the sums are constrained to
$N \geq j_{s+1}=j+\sum_{\ell =1}^s k_\ell > 0$ for $s<n$ (but it can be $j_0=j+\sum_{\ell =1}^n k_\ell = 0$).

\end{itemize}
When convergent, the sum in the r.h.s.\,of \eqref{dysexp} defines a solution of \eqref{eqhieraerror}
with initial condition $E_j(s)$.

For any sequence $k_1\dots k_n=\uk$, we denote by $n^+_{\uk}=n^+$ (resp.\;$n^-,n^=$) the number of times where $1$ (resp.\;$-1,-2$) appears in $\{k_1,\dots,k_n\}$, namely $n^+_{\uk}=\suml_{1\leq i\leq n}\delta_{1,k_i}$.
It is the number of operators $D^1_j$ (resp.\;$D^{-1}_j$, $D^{-2}_j$) appearing in the string $U_{j_1}(t,t_2)D_{j_1}^{k_1} \cdots D_{j_n}^{k_n}U_{j_0}(t_n,s)$. We have that
\be\label{jjo}
j_0=j+n^+-n^--2n^==j+n^+-m^-, \qquad n=n^++n^-+n^= 
\ee
where $m^-=n^-+2n^=$ is the number of negative steps performed by the process.

Note that the r.h.s.\;of \eqref{dysexp} is a (finite dimensional) functional integral over the space of all random walks with single positive and single and double negative jumps.
where we assume that $k_1, \cdots, k_n$ satisfy the following constraint.
For all integers $s$ s.t.\;$1 \leq s <n$,

$$j+\sum_{i=1}^s  k_i >0\;; \,\,\,\,\, j+\sum_{i=1}^n  k_i  \geq 0\;,$$
and both quantities are not larger than $N$.

\begin{Prop}\label{prop1}

 Suppose that,
for  some $A_0\geq 1$,
\be\label{A}
\| E_j(s)\|_1 \leq \left(\frac {j^2}N\right)^{j/2} A_0^j,\ \mbox{ for all }j
=1,\dots,N.
\ee
Then there exists $\t_0$ sufficiently small such that for any $\t=t-s \leq \t_0$, 
\be\label{A1}
\| E_j(t) \|_1 \leq \left(\frac {j^2}N\right)^{j/2} (A_1)^j,\ \mbox{ for all } j=1,\dots,N,
\ee
where $A_1=C(\t_0) A_0$, for some explicitly computable constant $C(\t_0) \geq 1$, depending only on $\t_0$.
\end{Prop}
\noindent
This Proposition is the heart of the paper. The proof of Theorem \ref{main} at the end of the present section   will consist in iterating this result, together with a scaling argument in order to remove the simplifying assumption \eqref{normV}.

\begin{proof}  The strategy of the proof of Proposition \ref{prop1} will be to split the sum over $k_1,k_2,\cdots, k_n$ in the Dyson expansion \eqref{dysexp} in several parts, and to use alternatively the three estimates at disposal 
\begin{itemize}
\item
$\| E_{j} \|_1 \leq 
\mbox{card}\{K\subset J=\{1,\dots,j\}\}
=2^{j}\;;
$
\item 
$\norm{D_j^{-1,-2}}\leq \frac{j^2}N$\;;
\item
$\norm{D_j^{1,-1,-2}}\leq  j$\;.
\end{itemize}

We have, for $j_0\neq 0$,
\bea\label{defdys}
&\| U_{j_1}(t,t_1)D_{j_1}^{k_1}U_{j_2} (t_1,t_2) \dots D_{j_n}^{k_n} U_{j_0}(t_n,s) \left(E_{j_0} (s)\right) \|_1  \nn \\
&\leq
\prod_{i=0}^{n} e^{ (j+k_i+\cdots + k_1) (t_{i}-t_{i+1}) }\prod_{i=1}^{n}  \| D_{j_i}^{k_i} \|  \|E_{j_0} \|_1\nn\\
&\leq  \prod_{i=0}^{n} e^{ (j+i) (t_{i}-t_{i+1}) }\prod_{i=1}^{n}  \| D_{j_i}^{k_i} \|  \|E_{j_0} \|_1
\eea 
($t_0=t$ and $t_{n+1}=s$).
For $j_0=0$, as already mentioned, one has to replace in \eqref{dysexp} the corresponding quantities defined in \eqref {E0}. 

Using that
\be\label{readily}
\prod_{i=0}^{n} e^{ (j+i) (t_{i}-t_{i+1}) }  \leq e^{j(t-s)} e^{n (t-s)}.
\ee
and $\int_s^{t=s+\tau} dt_1 \int_s^{t_1} dt_2 \cdots \int_s^{t_{n-1} }dt_n=\frac{\tau^n}{n!},$
\eqref{defdys} leads to
\be\label{dystot}
\| E_j (t) \|_1 \leq  e^{j\t} \sum_{n  \geq 0} \sum_{k_1 \cdots k_n} 
 \prod_{i=1}^{n}  \| D_{j_i}^{k_i} \|   \| E_{j_0} \|_1 \frac {\t^n}{n!} e^{n \t}.
\ee

We split now
\be
\sum_{k_1 \cdots k_n}=\sum_{k_1 \cdots k_n} \chi( j_0 <j) +\sum_{k_1 \cdots k_n} \chi (j_0 \geq j)\;,
\ee
namely, using \eqref{jjo},
\be
\sum_{k_1 \cdots k_n}=\sumls{k_1 \cdots k_n\\n^+-n^--2n^=<0}
  +\sumls{k_1 \cdots k_n\\n^+-n^--2n^=\geq 0}:=\sum_{k_1 \cdots k_n}^<
  +\sum_{k_1 \cdots k_n}^\geq.
\ee
The corresponding terms in the 
r.h.s.\,of \eqref{dystot} will be denoted by $\ej^<$ and $\ej^{\geq}$  respectively.

Let us first bound 
$$
\ej^<=e^{j\t} \suml_{n >0} \sumls{k_1 \cdots k_n}^<  \prodl_{i=1}^{n}  \| D_{j_i}^{k_i} \|   \| E_{j_0} \|_1 \frac {\t^n}{n!} e^{n \t}.
$$  
In this case $n^-+2n^=-n^+=m^--n^+=j-j_0>0$ { where $m^-$ is the number of negative jumps.} 

Therefore there must be $n_0$ such that 
 $\suml_{i=1}^{n_0}  k_i =0$,  $\suml_{i=n_0+1}^{s}  k_i  <0$ for all  $n_0+1 \leq s<n$
and  $\suml_{i=n_0+1}^{n}  k_i  = j_0 -j$.
This means that the random walk is definitively below
 $j$ from $n_0$ on.

The sequence $k_{n_0+1}\dots k_n$ has the associated numbers  $\tilde n^+,\tilde n^-,\tilde n^=$ satisfying
\be
\tilde n^++\tilde n^-+\tilde n^==n-n_0\ \ \ \mbox{ and }\ \ \ \tilde n^+-\tilde n^--2\tilde n^==j_0-j\;.
\ee
This implies that
\be
\tilde n^-+\tilde n^=\geq\frac{j-j_0}2\mbox{ and }
\tilde n^+\leq n-n_0-\frac{j-j_0}2.
\ee
Therefore, using \eqref{normd}, we get


 $$
 \|  D_{j_{n_0+1}}^{k_{n_0+1} }\|  \dots \| D_{j_n}^{k_n} \|  \leq (\frac {j^2} {N} )^{(j-j_0)/2}  j^{(n-n_0-(j-j_0)/2)}.
 $$
 Clearly the same estimate holds true also when $j_0=0$.
 
 On the other hand, observe that, for all $ \uk=k_1\dots k_n$ and $ i=1,\dots,n$,  $\| D_{j_i}^{k_i} \| \leq j_i
 \leq  (j+n^+)$, 
so that 
\be
\label{basic}
 \| D_{j_1}^{k_1} \| \dots \| D_{j_{n_0}}^{k_{n_0} } \| \leq   (j+n^+) ^{n_0}.
 \ee

 Thus we get, using \eqref{A}, and the obvious inequality
 $$
  j^{n-n_0} (j+n^+) ^{n_0} \leq (j+n^+) ^n,
 $$
\be
 \ej^<  \leq  e^{j\t} \sum_{n >0} \sum_{k_1 \cdots k_n}   \frac {\t_1^n}{n!}  (j+n^+) ^{n}  \left(\frac { j_0^2} N \right)^{j_0/2}  \left(\frac {j} {N} \right)^{ \frac {(j-j_0)}2}A_0^{j_0} ,
\ee
where 
$
\t_1=\t e^{\t}.
$

By using that  $\frac{k^k}{k!}\leq \sum\limits_{r \geq 0}\frac{k^r}{r!}= e^k $, we get
\be\label{stirinq}
\frac {\t_1^n}{n!} (j+n^+)^n \leq (  e\t_1)^n (\frac {j+n} n )^n \leq   ( e\t_1)^n e^j\;.
\ee
Note that $e\t_1=e\t e^\t$ is arbitrarily small provided that $\t_0\geq \t=t-s $ is sufficiently small.

In conclusion, since $A_0>1$ and $(\suml_{k_1 \cdots k_n}1)  \leq 3^n$, setting $\l =\l (\t_0)=3 e\t_0 e^{\t_0}$, we have
\bea\label{e<}
\ej^<  &\leq&  {e^{j (\t +1)}} A_0^j \sum_{n >0} \sum_{k_1 \cdots k_n} ( e\t e^\t)^n \left(\frac { j_0^2} N \right)^{\frac {j_0}2}  \left(\frac {j^2} {N} \right)^{ \frac {(j-j_0)}2}\nn\\
& 	\leq&
 (A_0 e^{\t +1})^j\frac { \l }{1-\l}
\left(\frac {j^2} N \right)^{j/2}.
\eea

As a consequence, if
$$
\lambda = 3 e\t_0 e^{\t_0 } \leq \frac 14,
$$  
we have
\be\label{condtau}
\frac { \l }{1-\l} \leq 1
\ee
and 
\be\label{t<}
\ej^<  \leq 
 \left(A_0 e^{\t +1}\right)^j 
\left(\frac {j^2} N \right)^{j/2}.
\ee

Next we estimate 
\be
\label{1}
\ej^\geq =e^{j\t} \sum_{n\geq 0}\sumls{k_1 \cdots k_n\\n^+-n^--2n^=\geq 0}
e^{n\tau}\frac {\t^n}{n!}   \prod_{i=1}^{n}  \| D_{j_i}^{k_i } \|  \| E_{j_0} \|_1.
\ee

\label{fofn}
To do this, we introduce a parameter $\mu \in (0,1) $, to be determined later,
and we split the above sum into two terms, namely 
$$
T_1=
 e^{j\t}\sum_{n \geq 0} \sumls{k_1 \cdots k_n\\n^+-n^--2n^=\geq 0\\ j_0   \leq { \mu \sqrt{N}}} e^{n\tau}\frac {\t^n}{n!}   \prod_{i=1}^{n}  \| D_{j_i}^{k_i} \|  \| E_{j_0} \|_1
$$
and
$$
T_2= e^{j\t} \sum_{n \geq 0} 
\sumls{k_1 \cdots k_n\\n^+-n^--2n^=\geq 0\\j_0 > { \mu \sqrt{N}}}
 e^{n\tau}\frac {\t^n}{n!}   \prod_{i=1}^{n}  \| D_{j_i}^{k_i} \|  \| E_{j_0} \|_1  \;.   
$$
Note that, when $j>{ \mu \sqrt{N}} $, $T_1=0$.

By \eqref{basic}, the inductive hypothesis \eqref{A} and estimate \eqref{stirinq} we deduce
\bea
T_1 &\leq&  e^{j(\t+1)} \sum_{n \geq 0}  \sum_{n \geq 0} \sumls{k_1 \cdots k_n\\n^+-n^--2n^=\geq 0\\ j_0   \leq { \mu \sqrt{N}}} ( \t e^{\t+1})^n \left( \frac { (j+\ell)^2} N\right)^{\frac {j+\ell}2} A_0^{(j+\ell)}\nn
\eea
where $\ell=n^+-m^-=j_0-j$. 
But
\bea
\big( \frac {(j+\ell)^2} N\big)^{\frac {j+\ell}2}A_0^\ell & = & \left( \frac {j^2}N \right)^{j/2} \left( \frac { (j+\ell)^2} N\right)^{\frac {\ell}2} \left( \frac { j+\ell } {j}\right)^jA_0^\ell \nn\\
& \leq & \left( \frac {j^2}N \right)^{j/2} \mu^\ell e^\ell A_0^\ell.
\eea
Here we used that $j+\ell =j_0 \leq \mu \sqrt {N}$ and that  $(\frac { j+\ell } {j})^j \leq e^\ell$.
Therefore we conclude that
\be\label{T1}
T_1 \leq   \frac 1 {1-\l} (A_0e^{\t+1})^j  \left( \frac {j^2}N \right)^{j/2} , 
\ee
after having fixed  $\mu=\frac 1 {A_0e}$.

For the term $T_2$ we make use of the estimate
\be\label{card}
\| E_{j} \|_1 \leq 
\mbox{card}\{K\subset J=\{1,\dots,j\}\}
=2^{j}
\ee
which is valid for any $j=1, \cdots, N$.
 Now we can assume $j< \frac {\mu \sqrt {N}} {2} $ because otherwise, reminding that $\mu= \frac 1 {eA_0}$,
 \be
 \label{bigj}
 (4A_0e)^j (\frac {j}{\sqrt {N}})^j >  (4A_0e)^j  (\frac {\mu}2)^j = 2^j \geq \| E_j \|_1.
 \ee
 Thus, if $j< \frac {\mu \sqrt {N}}2 $ and $j_0 > \mu \sqrt{N}$ then
 $$
 \mu \sqrt{N} \leq j_0=j+\suml_{i=1}^n k_i\leq j+n < \frac {\mu \sqrt {N}}2+n,
 $$  
 therefore
 $$
 n \geq \frac {\mu \sqrt{N}}2.
 $$ 
 As a consequence
 \be
T_2 \leq  2^j e^{j(\t +1)} \suml_{n \geq \frac {\mu \sqrt{N}}2}  (2\l)^n \leq  \frac {(2e)^{j(\t +1)}} {1-2\l} (2\l)^{\frac {\mu \sqrt {N}}2}. \nn
\ee
Setting $\b=|\log (2\l) |$ we obtain ($ \l \leq \frac 14$)
$$
(2\l)^{\frac {\mu \sqrt {N}}2} \leq (\frac 1{N})^{j /2} \sup e^{- \b \frac {\mu \sqrt {N}}2}  N^{j/2} \leq (\frac j{\sqrt{N}})^j (\frac 2{\b \mu})^{j }
$$
so that
\be
\label{eT2}
T_2 \leq 2(\frac {j^2} {N})^{j /2} (2e)^{j(\t +1)} (\frac 2{\b \mu})^{j }.
\ee

Collecting \eqref{t<},  \eqref{T1}, \eqref{bigj} and \eqref{eT2}, we conclude that there exists
an explicitly computable  
constant $C(\tau_0)\geq 1$
so that,  taking $\tau_0$ small enough
 the Proposition holds with $A_1=C(\tau_0)A_0$.   
\end{proof} 

To prove Theorem \ref{main}, we first fix $\tau_0$ small enough. Then, for $t<\tau_0$, Proposition \ref{prop1} gives the desired bound with
$C_1=0,\ C_2=C(\tau_0)C_0$. For $t\geq \tau_0$, we set
 $t=k\t$ with $\frac {\t_0} 2 \leq \t \leq \t_0$ and 
 $k\geq 1$ is a positive integer. Then we iterate  Proposition \ref{prop1} to get
$$
 \norm{E_j(t)
 }_1 \leq \left(\frac{j^2} N\right)^{j/2} (C(\t_0)^k C_0)^j \;.
$$
The first part of Theorem \ref{main}, under hypothesis \eqref{normV}, follows by setting  $C_1=4\frac {\log C(\tau_0)}{\t_0} $ and  $C_2=C(\tau_0)C_0$. 

In order to evaluate $C_1,C_2$ we first fix $\lambda(\tau_0)=3\tau_0ee^{\tau_0}=\frac14$. This easily implies that $|\log(1-\frac1{24e})|\leq\tau_0\leq 
\frac1{12e}(1-\frac1{24e})
$ . Moreover
tracing the dependence in $\tau$ of  \eqref{t<},  \eqref{T1}, 
 \eqref{bigj} and \eqref{eT2}, we easily show that, when $\lambda(\tau_0)=\frac14$, $C(\tau_0)$ can be taken as $C(\tau_0)=16(2e)^{\tau_0+1}\leq 16(2e)^{\frac1{12e}+1}$ so that one can take 
\be\label{c1c2}
C_1= 4\frac{16(2e)^{\frac1{12e}+1}}{|\log(1-\frac1{24e})|}
 \mbox{ and } 
  C_2={16(2e)^{1+1/12e}}{} C_0 
  \;.
\ee

What is left in order to finish the proof of  \bcrr{the first part of }Theorem \ref{main}, is to release the hypothesis \eqref{normV}. This is easily done by means of a scaling argument.  Note that rescaling $V$ as $V/\left(4\|V\|\right)$ in \eqref{eqhiera} and in \eqref{mfe} is equivalent to speed up time by $4\|V\|$ and rescale $K_0$ as $K_0/\left(4\|V\|\right)$.
But all the estimates is this section depend only on $K_0$ through the hypothesis \eqref{unitu}  
and we supposed that $\norm{e^{sK_0^j}}=\norm{e^{t\left(K_0^j+\frac {T_j}{N}\right)}}=1$ for all $s,t\in\bR$.
This allows to conclude.

In particular, for any $K_0, V$, the corresponding correlation error is 
$$E_j^{\left(K_0,V\right)}(t) = E_j^{\left(\frac{K_0}{4\|V\|},\frac{V}{4\|V\|}\right)}(4t\|V\|)$$ 
and we obtain
\be
\| E_j^{\left(K_0,V\right)}(t) \|_1 \leq \left( C_2 e^{C_1 t\norm{V}} \right)^j \left(\frac{j}{\sqrt N}\right)^j \;.
\ee

The second statement in Theorem \ref{main} is proven directly from the equation \eqref{eqhieraerror} which reads for $j=1$
\bea
\pa_t E_1&=& \left( K_0 +D_1\right) \left(E_1\right) 
 + 
D_1^1
\left(E_{2}\right) + 
D_1^{-1}\left(
E_{0}\right)\;.
\eea
 Using the semigroup $U_1(t,0)$, estimated by \eqref{gronwald} as $\norm{U_1(t,s)} \leq ~e^{(t-s)}$, 
 reminding that, under assumption \eqref{normV},
 $$ 
 \| D_1^{-1} \left(E_{0}\right)(s) \|_1 = \frac 1N \| Q(F,F) (s) \|_1 \leq  \frac{1}{4N}
 $$
and, by Theorem \ref{main},
$$
\norm{D_1^1 \left(E_{2}\right)(s)}_1  \leq  \norm {E_{2}(s)}_1 \leq 
 \frac{4}{N}C^2_2 e^{2C_1s/4}\;, 
$$
we get 
\bea
\norm{E_1(t)
}_1
&\leq&  e^t \| E_1(0) \|_1+
\int_0^t e^{(t-s)} \norm{D_1^1\left(
E_{2}\right)(s) + 
D_1^{-1}
\left(E_{0}\right)(s)}_1ds\nn\\
&\leq& e^t \frac{B_0}{N}+
\int_0^t e^{t-s} \left( \frac{1}{4N} + \frac{4}{N}C^2_2 e^{2C_1s/4} \right)ds\leq \frac{B_2e^{B_1t/4}}N
\nn\\
\mbox{with}& &B_1=2{C_1} \  \mbox{ and }\ B_2=B_0+\frac12+\frac{16C_2^2}{{C_1}-2}.
\label{b1b2} 
\eea
The same argument as before ($t\to 4t\norm{V}$) allows to place the $\norm{V}$ in the exponential, so
Theorem \ref{main} is proven.

\qed
\begin{Rmk}\label{srtucture}
Let us point out that the proof of Proposition \ref{prop1} out of the equation \eqref{eqhieraerror} uses  only  the  properties \eqref{unitu}, \eqref{normdj},
 \eqref{normd}
 and the ``convention" \eqref{E0}. No particular explicit incarnation for the operators $D_j$s in \eqref{eqhieraerror} is needed. Therefore  Theorem \ref{main} and Corollary \ref{cormain} remain true under the validity of Proposition \ref{maineq} involving an equation of type (\ref{eqhieraerror}-\ref{E0})  endowed with the following assumptions
$$
\norm{e^{\left(K_0^j+\frac {T_j}{N}\right)t}}=1,\ \| D_j \| , \| D^1_j\| \leq j  \mbox{ and }  \| D_j^{-1}\|,
\| D_j^{-2}\|, \|D^{-1}_1(E_0)\|,  \|D^{-2}_2(E_0)\|
 \leq \frac {j^2} N.
$$
\end{Rmk}

\begin{Rmk} According to the heuristic argument in the Introduction, we believe that our estimates are optimal as regards the size of chaos and the rate of convergence \bcb (see also the classical estimate on independent random variables for which the same result is easily obtained (e.g. \cite{Gr71,Sz89})).\ec They are certainly not optimal as regards the time dependence.
However in the above proof we did not try to optimize the numerical constants. It is easy to realize that such constants affect only the growth of the error as a function of time, but not the dependence on $j$ and $N$ whose analysis is the main purpose of this paper.
\end{Rmk}

\section{Return to the concrete examples} \label{applications} 

%
\bcb In this section we turn back to the concrete models we have in mind as expressed  by the table in Section \ref{ss:EE}. 
\ec

\subsection{Stochastic models}\label{km}

We recall the evolution equation for the probability measure $f^N (\calv_N,t)$ describing the  Kac model:
\bcb
\be\label{kacnbody}
\pa_t f^N(\cal V_N,t)=\frac 1N \sum_{i<j} \int d\omega B(\omega; v_i-v_j)\{f^N(\calv_N^{i,j},t)-f^N(\calv_N,t)\}\;,
\ee
\ec
where $\calv_N=\{v_1, \cdots ,v_N\}$ and  $\calv_N^{i,j}=\{v_1, \cdots,v_{i-1}, v_i', \cdots, v_{j-1},v_j', \cdots, v_N\}$  is the vector of the velocities after a collision between particle $i$ and $j$ and
$\tfrac{B(\omega; v_i-v_j)}{|v_i-v_j|}$ is the differential cross-section of the two-body process, which we assume here to be bounded.
The resulting kinetic equation reads
\bcb
\be\label{vlasovkac}
\pa_t f(v,t) = \int dv_1 \int d\omega B(\omega; v-v_1) \{ f(v',t) f^t(v_1',t) -f(v,t) f^t(v_1,t) \}\;.
\ee
\ec
 
 For the soft spheres model
the  probability density $f^N(X_N, \calv_N,t)$ evolves according to
\bcb
\bea\label{softspheresNbody}
\pa_t f^N + \suml_{i=1}^N v_i \cdot \nabla_{x_i} f^N &=&
\frac 1N \suml_{i<j} \, h\left(|x_i-x_j|\right) B\left(\frac{x_i-x_j}{|x_i-x_j|}; v_i-v_j\right)
\nn \\
&&
\times\{f^N(X_N ,\calv_N^{i,j},t)-f^N(X_N,\calv_N,t)\}\;. 
\eea
\ec
Here $X_N= \{x_1, \cdots, x_N \}$ and $h:\R^+ \to \R^+$ is a positive function with compact support. 
The associated kinetic equation is \bcb
\bea\label{vlasovsoftsphesres}
\pa_t f(x,v,t) + v\cdot \nabla_{x}f (x,v,t) &=&\int dv_1 \int dx_1\, h(|x-x_1|) B\left(\frac{x-x_1}{|x-x_1|}; v-v_1\right)\nn   \\
&\times &\{ f^t(x, v',t) f^t(x_1, v_1',t) -f^t(x, v,t) f^t(x_1, v_1,t) \}.  
\eea
\ec

For both models the correlation error is defined as
\bcb
\be\label{corkac}
E_j (t) := \sum_{K \subset J} (-1)^{|K|}  f^N_{J \backslash K} (t) f(t)^{\otimes K}
\ee
\ec
where $J=\{ 1,2, \cdots, j \}$, $K$ is any subset of $J$ and $|K| =$ cardinality of $K $.
$f_A ^N(t)$ stands for the $|A|-$marginal $f_{|A|} ^N(t)$ computed in the configuration $ \{z_i \}_{i\in A} $. 
Similarly, $f(t)^{\otimes K} = f(t)^{\otimes |K|}$ evaluated in $ \{z_i \}_{i\in K}$.
Moreover either $z_i=v_i$ or  $z_i=(x_i,v_i) $ for the Kac and soft spheres model respectively
and  $f$ is the solution to the kinetic equations  above. 


According to Corollary  \ref{cormain} we have \bcb the following result.
\begin{Thm}\label{smain}
Let us suppose that $B$ is bounded and that the initial conditions of equations \eqref{kacnbody} and \eqref{softspheresNbody} lead, through \eqref{corkac} {for some probability density $f(0)$}, to quantities $E_j(0),j=2,\dots,N$ satisfying \eqref{hyp} together with $\|E_1(0) \|_{L^1 (dz)} \leq B_0/N$.

Then, 
for both Kac and soft sphere models and 
 for all $t\in\bR \mbox{ and }j=\unN$, the marginals of the solution $f^N(t)$ of \eqref{kacnbody} or \eqref{softspheresNbody} 
 satisfies
\be
\label{stoccor}
\| f^N_j (\cdot,t) -f^{\otimes j}  (\cdot,t) \|_{L^1 
} \leq 
D_2 e^{
2
D_1t 
\| B \|_{\infty}\left(2\p + \| h \|_{\infty}\right)
} \frac {j^2}N\;,
\ee
where $f(t)$ is the solution of equations \eqref{vlasovkac} or \eqref{vlasovsoftsphesres} respectively, with initial condition $f(0)$. 

Here 
 $D_1,D_2$ are the geometrical constants defined  in Corollary \ref{cormain}. 
\end{Thm}
\ec
%

The Kac model has been extensively studied. We refer to the recent papers  \cite{MM,MMW,HaMi} and 
 to the references therein.
Here uniform in time estimates have been derived for models including unbounded kernels.
Typically the error in the propagation of chaos is controlled in terms of a Wasserstein distance, which however is
sensitive on the  metric chosen for the configuration space. 

Estimates in $L^1$ similar to \eqref{stoccor} for  models with bounded cross-section 
were obtained first in \cite{GM97}. 
The technique uses an explicit representation of the underlying stochastic process, making rigorous the 
heuristic argument described in our introduction. 

In contrast, the method of the present paper focuses on the errors $E_j$ in an abstract setting, thus using only the hierarchy  of equations. 
This allows us to apply our results also to different cases such as the soft sphere model
where the impact parameter is not random.
We also remark that collisional mean field models are potentially useful for applications in population dynamics involving a large number of agents.

\subsection{Quantum mean field}\label{final}
In this case, $\AA= \mathcal L^1(L^2(\bR^d))$ is the space of trace-class operators on $L^2(\bR^d)$ and
\begin{itemize}
\item $K=\frac1{i\hbar}[-\frac{\hbar^2}{2}\Delta_{\bR^d},\cdot]$, $K^N=\frac1{i\hbar}[-\frac{\hbar^2}{2}\Delta_{\bR^{Nd}},\cdot]\;;$
\item $V_{r,j}=\frac1{i\hbar}[V(x_r-x_j),\cdot]$, 
$V^N=\frac1{i\hbar}[ \frac{1}{N}\suml_{1\leq i\neq j\leq N}V(x_i-x_j),\cdot]\;.$
\end{itemize}
For any operator  $B\in \AA_j $, $j=\unN$, we denote its integral kernel by $B(X_j,X_j')=~B(Z_j),\ X_j=~(x_1,\dots,x_j)\in\bR^{jd}$ where we  denote $Z_j=(z_1,\dots,z_j)\in\bR^{2jd},\ z_k=(x_k,x'_k)\in\bR^{2d}$ . We also denote  $Z_{J \backslash K}$ the vector $Z_j$ after removing the components $z_{i_1},\dots z_{i_k}$, where $J$ stands for $\{1,\cdots,j\}$ and $K = \{i_1,\cdots,i_k\}$.
The  formulas defining the error and its inverse read
\bea
E_j(Z_j)&=&
\sum_{k=0}^j\sum_{1\leq i_1,\cdots,i_k\leq j} \frac{1}{k!}(-1)^k
\FF(z_{i_1})\dots F(z_{i_k})F^N_{j-k}(Z_{J \backslash K})\;,\nn\\
\FF^N_j(Z_j)&=&\sum_{k=0}^j\sum_{1\leq i_1,\dots,i_k\leq j}\frac{1}{k!}
\FF(z_{i_1})\dots F(z_{i_k})E_{j-k}(Z_{J \backslash K})\;,\nn
\eea 
\bcb
where the marginals $F^N_{j}$ are defined through their integral kernel
\be
F^N_{j}(Z_j)=\int_{\bR^{(j-k)d}} F^N(Z_j,\bar Z_{j-k})d\bar Z_{j-k}\nn
\ee
where $\bar Z_{j-k}=((x_{j+1},x_{j+1}),\dots,(x_{N},x_{N})\in\bR^{(j-k)d}$.
\ec

\noindent Note that all operators commute since they act on different variables.

The quantum $N$-body dynamics is defined by the equation
\bea\label{nbodyqaunt}
i\hbar\partial_t \FF^N(t)&=&[-\frac{\hbar^2}{2}\Delta_{\bR^{Nd}}+ \frac{1}{N}\sum_{1\leq i<j\leq N}V(x_i-x_j),\FF^N(t)]\\
\FF^N(0)&\geq & 0,\  \tr F^N(0)=1.\nn
\eea
Moreover the mean field Hartree equation reads
\be\label{hartree}
\partial_tF=\frac1{i\hbar}[-\frac{\hbar^2}{2}\Delta+V_F,F],
\ee
where 
\be\nn
V_F(x)=\int_{\bR^d}V(x-y)F(y,y)dy.
\ee

%

Corollary \ref{cormain} is reformulated as follows \bcb (we state here only the case on factorized initial data. The reader interested can easily extend the result to the general case).
\begin{Thm}\label{qmain}
Let us suppose that $V$ is bounded and that the initial condition $F^N(0)$ of the N-body quantum problem \eqref{nbodyqaunt} satisfies
\be\nn
F^N(0)=F(0)^{\otimes N},\ F(0)\geq 0,\ \tr{F(0)}=1,
\ee
and let $F(t)$ the solution of \eqref{hartree} with initial condition $F(0)$.

Then, for all $t\in\bR, 0<\hbar\leq 1\ \mbox{(say), and }j=\unN$, 
\be\label{cormainquant}
\tr|\FF^N_j (t) - \FF^{\otimes j}  (t)| \leq 
D_2 e^{2D_1t\frac{\norm{V}_\infty}\hbar}
\ \frac {j^2}{{ N}}.
\ee
Here the constants $D_1,D_2$ are the ones of Corollary \ref{cormain}.
\end{Thm}

A bound like  (\ref{cormainquant}) is new to our knowledge, \ec in view of its explicit dependence in 
$j$, {$t$ and $\hbar$ } and validity for ``all" (pure and mixed) initial data.

Introduced in 1927 (one year after the Schr\"odinger equation), the Hartree equation has received an 
enormous interest in physics since then. The first derivation from the quantum $N$-body dynamics of observables goes back to Hepp in \cite{hepp}, 
using coherent states, and to Spohn in \cite{Spohn1} for pure states, using hierarchies. 
A proof of the mean field limit for bounded potentials and mixed states including rates of convergence can be found in \cite{BGM}, 
and for Coulomb singularity and pure states in \cite{EY}, after \cite{BEGMY}, using hierarchies. 
At $\hbar=1$, the rate of convergence in $\frac1{\sqrt N}$ has been discussed in \cite{kpickl}
together with an explicit dependence in $j$, using heavily a  pure states hypothesis on the initial data,
and in \cite{BPS}, improved to a (optimal) rate in $\frac1N$ in \cite{CLS} without tracing the dependence 
in $j$ and $\hbar$. These two papers use the method of second quantization and ``coherent" states 
in Fock space, initiated in \cite{hepp,GV}. A rate of convergence with explicit (exponential) 
dependence in $j$ and $\hbar$ can be found in \cite{GPP}, Theorem 7.1. Eq.\,\eqref{cormainquant} realizes an 
improvement of this last result.

\section*{Appendix. Derivation of the correlation equations}


We prove  here Proposition \ref{maineq}.

We want to compute the time-derivative of \eqref{deferror}, which we recall:
\be
E_j=\sum_{K \subset J} (-1)^{k} \ 
\FF^{\otimes K,J} \FF^J_{J \backslash K}
\ee
with $J=\{1,2 \dots j\}$, $k = |K|$.

We first notice that, using 
the mapping \eqref{o} and the hierarchy \eqref{eqhiera},
\bea
\label{eveq}
\FF^{\otimes K,J}\pa_t F^J_{J \backslash K} &=&\left(K_0^{J \backslash K} +  \frac {T_{J \backslash K}}{N}\right) \left(\FF^{\otimes K,J}F^J_{J \backslash K}\right)\nn\\
&& + \a(j-k,N) C_{J \backslash K, j+1} \left( \FF^{\otimes K,J
\cup \{j+1\}} F^{J\cup \{j+1\}}_{(J \backslash K) \cup \{j+1\}}\right)\;,
\eea
where, for $S \subset J$,
\be
K^S_0= \sum_{i \in S} \mathbb I_{ \mathfrak{A}_{i-1}}\otimes K_0\otimes\mathbb I_{\mathfrak{A}_{j-i}}\;, \,\,\,\,\,\, T_S=\sum_{\substack{i,r \in S \\ i < r}} T_{i,r}
\ee
and
\be
C_{S, j+1} =\sum_{i\in S} C_{i,j+1}\;.
\ee

Moreover, \eqref{mfe} and \eqref{ODt} imply
\be
\partial_t F^{\otimes K,J} = K_0^K \left(F^{\otimes K,J}\right) + \sum_{i \in K} F^{\otimes K^i,J}
Q(F,F)^{\otimes \{i\},J}\;,
\ee
with the notation $K^{i} = K \backslash\{i\}$.

Therefore we have
\bea
\label {eqE}
\pa_t E_j=& K_0^J \left(E_j\right) +\suml_{K \subset J} (-1)^{k} \suml_{i\in K}  F^{\otimes K^i,J} Q(F,F)^{\otimes\{i\},J} F^J_{J
\backslash K} \nn \\
& +  \suml_{K \subset J} (-1)^{k} \a(j-k,N) \suml_{i\in J \backslash K} C_{i,j+1} \left(F^{\otimes K,J\cup\{j+1\}} F^{J\cup \{j+1\}}_{(J \backslash K) \cup \{j+1\}}\right)  \nn \\
& +\frac 1{2N} \suml_{K \subset J} (-1)^{k} \suml_{\substack{i,r \in J\backslash K \\ i \neq r}} T_{i,r} \left(F^{\otimes K,J}F^{J}_{J \backslash K}\right)\;.
\eea
In the following computation we shall simplify the notation by skipping the upper indices of sets $J$ and $J\cup \{j+1\}$, now clear from the context.

We compute next the terms in the three lines on the r.h.s.\,of \eqref{eqE} separately. They are denoted by ${\cal T} _i$, $i=1,2,3$
respectively.

Using $K^i= K \backslash\{ i \}$ and the change $ K \to K^i$,
\bea
{\cal T} _1=&  K_0^J\left(E_j\right)+\suml_{i\in J} \suml_{K \subset J^{i}} (-1)^{(k-1)} F^{\otimes K} Q(F,F)^{\otimes\{i\}} 
F_{J^i \backslash K}  \nn \\
&=K_0^J\left(E_j\right)-\suml_{i\in J} Q(F,F)^{\otimes\{i\}}  \suml_{K \subset J^{i}} (-1)^{k} F^{\otimes K} 
F_{J^i \backslash K}  \nn \\
&=K_0^J\left(E_j\right) - \suml_{i\in J} Q(F,F)^{\otimes\{i\}}  E_{J^i}\;,\nn
\eea
where in the last step we applied \eqref{eq:dakaKb} (for $E^J_{J^{i}}$).

To compute the term ${\cal T}_2$, we will make use of the combinatorial identity
\be
\label{identity}
\sum_{ K \subset R} (-1)^k \a ( j-k,N) = \a (j,N) \delta_{R, \emptyset }-\frac 1N  \delta_{|R|,1}
\quad R \subset J\;.
\ee
We postpone the elementary proof of \eqref{identity} to the end of the section.

Applying again \eqref{eq:dakaKb} into the second line of \eqref{eqE}
we obtain
\bea
\label{T2}
{\cal T} _2 &=&
 \sum_{K \subset J} (-1)^{k} \a(j-k,N) \sum_{i\in J/K} C_{i,j+1} 
\left( F^{\otimes K}\sum_{L \subset J \cup\{j+1\}\backslash K } F^{\otimes  L}
E_{J\cup \{ j+1 \} \backslash (K\cup L) } \right)\nn \\
&=&\sum_{i\in J}  \sum_{K \subset J^i} (-1)^{k} \a(j-k,N) C_{i,j+1} \left( F^{\otimes K}\sum_{L \subset J \cup\{j+1\}\backslash K } F^{\otimes  L}
E_{J\cup \{ j+1 \} \backslash (K\cup L) } \right)\;.\nn \\
\eea
Now we distinguish the following cases:
\begin {itemize}
\item  $r=1$: \quad  $i, j+1 \in L$
\item   $r=2$:  \quad $i, j+1 \notin L$
\item   $r=3$:  \quad $i \in L , j+1 \notin L$
\item   $r=4$:  \quad $i \notin L , j+1 \in L$
\end{itemize}
and set
$$
{\cal T}^1 _2=\sum_{r=1}^4 {\cal T}^r _2
$$
with ${\cal T}^i _2$ given by \eqref{T2} with the additional constraint $r=i$ above. Setting
$L'=L^{i,j+1} = L \backslash \{i,j+1\}$ and $R=K \cup L'$ and recalling \eqref{eq:EI},
\bea
{\cal T}^1 _2
 &=  \suml_{i\in J} \suml_{R \subset J^i }Q(F,F)^{\otimes\{i \}}  F^{\otimes R}   E_{J^{i}\backslash R} \suml_{K \subset R} (-1)^{k} \a (j-k,N)\nn
 \eea
 so that  \eqref{identity} leads to
 \be
 {\cal T}^1 _2 =\a (j,N)  \sum_{i\in J}  Q(F,F)^{\otimes\{i \}} E_{J^i}- 
 \frac 1N   \suml_{\substack{i,r \in J \\ i \neq r}} Q(F,F)^{\otimes\{i \}}F^{\otimes \{r\}} E_{J^{i,r}}\;.\nn
\ee

Observe that
$$
{\cal T}_1+{\cal T}^1 _2 =K_0^J \left(E_j\right)-\frac jN  \sum_{i\in J} Q(F,F)^{\otimes\{i \}} E_{J^i}
- 
 \frac 1N   \sum_{\substack{i,r \in J \\ i \neq r}} Q(F,F)^{\otimes\{i \}}F^{\otimes \{r\}} E_{J^{i,r}}\;,
$$
namely there is a crucial compensation for which all the operators, but  $K_0^J $, are $O(\frac{j^2}{N})$.

Furthermore, setting $R=K \cup L$,
\bea
{\cal T} ^2_2=&   \suml_{i\in J} \suml_{L \subset J^i  }\suml_{K \subset J^i/L} (-1)^{k} \a(j-k,N)C_{i,j+1}
\left( F^{\otimes(K\cup L)}   E_{J\cup \{ j+1 \} \backslash (K\cup L) } \right) \nn \\
& = \suml_{i\in J} \suml_{R \subset J^i  } \suml_{K \subset R}(-1)^k \a(j-k,N)  C_{i,j+1}
\left(F^{\otimes R} E_{ J \cup \{j+1\} \backslash R}\right) \nn \\
& = \a(j,N) \suml_{i\in J} C_{i,j+1}  \left(E_{ j+1 }\right) 
-\frac 1N \suml_{\substack{i,r \in J \\ i \neq r}}  C_{i,j+1}\left( F^{\otimes\{r\}}E_{ J^r \cup \{j+1 \}}\right)\;. \nn
\eea

To compute ${\cal T} ^3_2$ we set $L'=L^{i}$ and $R= L' \cup K$ so that

\bea
{\cal T}^3 _2
&=& \sum_{i\in J}  \sum_{R \subset J^i } \sum_{K \subset R } (-1)^{k} \a(j-k,N)C_{i,j+1}\left(F^{\otimes R} F^{\otimes \{i\}}   
E_{ J^i\cup \{ j+1 \} \backslash R}\right)\nn \\
&=& \a(j,N) \sum_{ i\in J}  C_{i,j+1} \left( F^{\otimes \{i\}} E_{J^i \cup \{ j+1\} } \right)
-\frac 1N  \sum_{\substack{i,r \in J \\ i \neq r}} C_{i,j+1} \left(F^{\otimes \{r,i\}} 
E_{ J^{i,r} \cup \{j+1 \}}\right)\;. \nn
\eea

Finally, setting $L'=L^{j+1}$ and $R=K \cup L'$ we obtain
\bea
{\cal T}^4 _2& =\suml_{i\in J} \suml_{R \subset J^i }  \suml_{K \subset R } (-1)^{k} \a(j-k,N)
C_{i,j+1} (F^{\otimes R}F^{\otimes \{j+1\}}  E_{J \backslash R }) \nn \\
& =\a(j,N) \suml_{ i\in J}  C_{i,j+1} ( F^{\otimes  \{j+1\}} E_{J}) -\frac 1N \suml_{\substack{i,r \in J \\ i \neq r}}  C_{i,j+1} 
\left(F^{\otimes\{r,j+1\}} 
E_{ J^{r} }\right)\nn  \;.
\eea

Similarly we compute the term
\bea
{\cal T} _3&=&\frac 1{2N} \sum_{K \subset J} (-1)^{k}  \sum_{ i,s \in J \backslash K} 
\sum_{L \subset J \backslash K}T_{i,s} \left(F^{\otimes K,J}F^{\otimes L}E_{J \backslash (K\cup L)}\right)\nn\\
&=& \sum_{r=1}^3 {\cal T}^r _3 \nn
\eea
where each term ${\cal T}^r _3$ corresponds to the constraints
\begin {itemize}
\item  $r=1$: \quad  $i, s \in L$
\item   $r=2$:  \quad $i, s \notin L$
\item   $r=3$:  \quad $i \in L , s \notin  L$\;.
\end{itemize}
Setting $L'=L^{i,s}$ and $K\cup L'=R$ we have
\bea
{\cal T}^1 _3=&\frac 1{2N}  \suml_{i,s \in J}  \suml_{R\subset J^{i,s}} 
\suml_{K \subset R} (-1)^{k} 
 T_{i,s} \left(F^{ \otimes\{i,s\}} F^{\otimes R}  E_{J^{i,s}\backslash R}\right)\nn \\
&=\frac 1{2N}  \suml_{i,s \in J}  T_{i,s} \left(F^{ \otimes\{i,s\}}  E_{J^{i,s}}\right) \nn
\eea
and an analogous computation gives
\be
{\cal T}^2 _3=\frac 1{2N}  \suml_{\substack{i,s \in J \\ i \neq s}} T_{i,s} \left(E_{j} \right) 
= \frac{T_j}{N}\left(E_j\right)\nn
\ee
and
\be
{\cal T}^3 _3= 
 \frac 1{N} \suml_{\substack{i,s \in J \\ i \neq s}}  T_{i,s} \left({\FF}^{\otimes\{i\}}E_{J^{i}}\right)\;.\nn
\ee

In conclusion:

\bea
\pa_t E_j=& K_0^j \left(E_j\right)+\frac {T_j}{N} \left(E_j\right) +\a(j,N) \suml_{ i\in J}  C_{i,j+1} \left( F^{\otimes\{ i\}}  E_{J^i \cup \{ j+1\} }+ F^{\otimes \{ j+1\}} E_{J}\right) \nn \\
&-\frac 1N   \suml_{\substack{i,s \in J \\ i \neq s}}  C_{i,j+1} \left(F^{\otimes \{ s\}}  E_{J^s\cup \{j+1 \} }\right) \nn \\
&+ \a(j,N) C_{j+1}\left( E_{j+1}\right) \nn \\
&-\frac jN  \suml_{i\in J} Q(F,F)^{\otimes\{ i\}}  E_{J^i}+ \frac 1N \suml_{\substack{i,s \in J \\ i \neq s}}
T_{i,s} \left(F^{\otimes\{ i\}}  E_{J^i} \right) \nn \\
&-\frac 1N \suml_{\substack{i,s \in J \\ i \neq s}}
 C_{i,j+1} \left(F^{\{ i,s\}}  E_{J^{i,s} \cup \{j+1\}}  +F^{\otimes\{ s,j+1\}}E_{J^s} \right)
\nn \\
& \frac 1{2N}\suml_{\substack{i,s \in J \\ i \neq s}} T_{i,s} \left(F^{\otimes\{ i,s\}}  E_{J^{i,s}}\right) - \frac 1{N} \suml_{\substack{i,s \in J \\ i \neq s}} Q(F,F)^{\otimes\{i\}} F^{\otimes\{ s\}}  E_{J^{i,s}}\;.
\label{eq:finApp}
\eea

We organized the terms in
the above equation, according to the following rule. The first two lines contain  operators which do not change the particle number. 
The third line increases the number of particles by one. The  $4^{th}$ and $5^{th}$ lines decrease the particle number by one. Finally the last line decreases it by two.

Using the definition of $D_j$ in \eqref{eqdefDj}-\eqref{eqdefDj'} and of $D_j^1, D_j^{-1}, D_{j}^{-2}$
in \eqref{eqdefDj1}, \eqref{eqdefDj-} and \eqref{defd=}, Eq.\,\eqref{eq:finApp} reads
\bea
\pa_t E_j&=& \left(K_0^j+\frac{T_j}{N}\right)\left(E_j\right) +
D_j
\left(E_{{j}}\right) \nn \\
&&+ 
D_j^1
\left(E_{j+1}\right) + 
D_j^{-1}
\left(E_{j-1}\right) +
D_j^{-2}
\left(E_{j-2}\right)\;.
\eea

We end this section with the proof of \eqref{identity}. Denoting $|R| = r$ one has
\bea
\sum_{k=0}^r (-1)^k \binom {r}{k} \a (j-k,N)&=& \a(j,N) \sum_{k=0}^r (-1)^k \binom {r}{k} +
\frac 1N \sum_{k=0}^r (-1)^k \binom {r}{k} k \nn \\
&=&\a(j,N) \delta_{r,0} - \frac 1N r \sum_{k=0}^{r-1} (-1)^k \binom {r-1}{k} \nn \\
&=& \a(j,N) \delta_{r,0} - \frac 1N r(1-1)^ {r-1}\nn\\
&=&\a(j,N) \delta_{r,0} - \frac 1N r\delta_{r,1}\;. \nn
\eea

\vskip 1cm
{\small
\textbf{Acknowledgements.} The authors are grateful to Joaquin Fontbona and St\'{e}phane Mischler for helpful discussions.
This work has been partially carried out thanks to the supports of the LIA AMU-CNRS-ECM-INdAM
Laboratoire Ypatie des Sciences Math\'ematiques (LYSM) and the A*MIDEX project (n$^o$ ANR-11-IDEX-0001-02) funded by the ``Investissements d'Avenir" French Government program, managed by the French 
National Research Agency (ANR). T.P. thanks also the Dipartimento di Matematica, Sapienza Universit\`a di Roma, for its kind hospitality during the elaboration of this work. S.\,S.\,acknowledges support of the German Research Foundation (DFG n$^o$ 269134396).
}

\end{document}